\documentclass[11pt]{article}

\usepackage{amsmath, amsfonts, amsthm, amssymb, mathtools}
\usepackage{authblk, color, bm, bbm, graphicx, epstopdf}
\usepackage[labelfont = bf, labelsep = period]{caption} %Option: font = small
\usepackage{subcaption, xspace, enumitem}
\usepackage{booktabs, tabularx, multirow, rotating, array}
\usepackage[bottom,hang]{footmisc}
\newcolumntype{C}{>{\centering\arraybackslash}X}
\newcolumntype{R}{>{\raggedleft\arraybackslash}X}
\newcolumntype{L}{>{\raggedright\arraybackslash}X}
\usepackage[flushleft]{threeparttable}
\usepackage{accents}
\DeclareRobustCommand{\ubar}[1]{\underaccent{\bar}{#1}}
\makeatletter
\newcommand{\mylabel}[2]{#2\def\@currentlabel{#2}\label{#1}}
\makeatother
\usepackage{tikz,pgfplots}
\usetikzlibrary{positioning,arrows.meta}
\pgfplotsset{
	every axis/.append style={
		height=4cm,
		width=8cm,
		ytick={3,9},
		yticklabels={$s=2$,$s=1$},
		xticklabels={,0,5,10,15,20,25},
		xmin=0,xmax=25,ymin=0,ymax=12,
		minor x tick num=4
	},
	axis line style={gray},
	every x tick label/.append style={font=\scriptsize}
}

\usepackage{fullpage}[2cm]
\title{A Scenario Decomposition Algorithm for Strategic Time Window Assignment Vehicle Routing Problems}

\author{Anirudh Subramanyam}
\author{Akang Wang}
\author{Chrysanthos E.~Gounaris}
\affil{\small Department of Chemical Engineering, Carnegie Mellon University, United States}
\date{\vspace{-0.1em}}

\begin{document}

\maketitle

\begin{abstract}
We study the strategic decision-making problem of assigning time windows to customers in the context of vehicle routing applications that are affected by operational uncertainty.
This problem, known as the Time Window Assignment Vehicle Routing Problem, can be viewed as a two-stage stochastic optimization problem, where time window assignments constitute first-stage decisions, vehicle routes adhering to the assigned time windows constitute second-stage decisions, and the objective is to minimize the expected routing costs.
To that end, we develop in this paper a new scenario decomposition algorithm to solve the sampled deterministic equivalent of this stochastic model.
From a modeling viewpoint, our approach can accommodate both continuous and discrete sets of feasible time window assignments as well as general scenario-based models of uncertainty for several routing-specific parameters, including customer demands and travel times, among others.
From an algorithmic viewpoint, our approach can be easily parallelized, can utilize any available vehicle routing solver as a black box, and can be readily modified as a heuristic for large-scale instances.
We perform a comprehensive computational study to demonstrate that our algorithm strongly outperforms all existing solution methods, as well as to quantify the trade-off between computational tractability and expected cost savings when considering a larger number of future scenarios during strategic time window assignment.

\end{abstract}

\section{Introduction}\label{sec:introduction}
The commitment to deliver (or pickup) goods within scheduled time windows is a common practice in several real world distribution networks. In many industries, these time windows are mutually agreed upon by the distributor and customer through long-term delivery contracts. For example, in a distribution network of retailers, it is common that deliveries to a retail store are always made on the same day of the week (at about the same time) for an entire year~\cite{Spliet2015_continuous,Vercammen2016}.
Likewise, in maritime distribution of liquefied natural gas, a central planning activity is to design and negotiate contractual agreements of \textit{annual delivery plans} that specify delivery dates and corresponding delivery quantities to customers~\cite{Zhang2015}.
From the customer's point of view, this is crucial for efficient inventory management and scheduling of personnel to process the delivery. From the distributor's point of view, it reduces the variability across repetitive deliveries and exposes efficiencies that add up to significant cost savings.
Short- and medium-term contracts of similar nature can be also found in small-package shipping where, for instance, courier companies provide a delivery time window to customers receiving sensitive packages~\cite{Jabali2015}.
Other examples of applications where such operations are typical include, among others, attended home delivery in e-commerce businesses~\cite{Agatz2011} and internet installation services~\cite{Ulmer2017}.

Once a time window has been agreed upon and communicated to the customer, the distributor must attempt to meet it on an operational (e.g., daily) basis as well as possible.
This is done by solving a Vehicle Routing Problem with Time Windows (VRPTW) to determine a delivery schedule that adheres to the agreed time windows. The assigned time windows strongly influence the structure of feasible delivery schedules and, hence, the daily incurred distribution costs.
Therefore, a natural choice is to assign time windows based on the arrival times at customer locations in the optimal (i.e., minimum cost) vehicle routing schedule. However, this seemingly optimal decision may become highly suboptimal in the presence of operational uncertainty.

In reality, operational level information (such as customer demands or travel times) is often not known with certainty at the strategic level when time windows are to be decided. For example, the demand volume of a customer typically fluctuates per delivery. Similarly, travel times vary on a day-to-day basis (e.g., because of unpredictable traffic conditions). The true values of these operational parameters are not known far in advance, and often may become known only on the day of delivery before the vehicles are dispatched. This makes the strategic assignment of time windows a non-trivial task.
Indeed, if one utilizes only nominal values of the uncertainty when assigning time windows, then it will often lead to situations in which the distribution costs are unacceptably high, since the nominal delivery schedule may no longer be feasible, let alone optimal, in such cases.
Fortunately, with the increasing availability of data, distributors can readily obtain forecasts of uncertain operational parameters (e.g., as perturbations from their nominal values). It is possible to take advantage of this information and assign time windows in a way that will lead to low distribution costs in the long run.
The goal of this paper is to study the problem of strategic time window assignment in the presence of operational uncertainty.

Our paper builds upon the work of~\cite{Spliet2015_continuous}, which introduced the Time Window Assignment Vehicle Routing Problem (TWAVRP).
The TWAVRP consists of assigning time windows of pre-specified width within some exogenous time windows to a set of known customers.
The exogenous time windows typically correspond to operating hours of the customer but may also arise from hours-of-work or other government regulations.
The work of~\cite{Spliet2015_continuous} studies the TWAVRP under situations in which the demand volume of the customers is unknown and subject to uncertainty. However, a finite set of ``scenarios,'' each describing a possible realization of demand for every customer, is assumed to be given with known probability of occurrence.
This information is used to formulate a two-stage stochastic program, in which the first-stage decisions are to assign time windows, while the second-stage decisions are to design vehicle routing schedules satisfying the assigned time windows, one for each of the demand scenarios. The objective is to minimize the total routing costs, averaged over the postulated scenarios.
A similar modeling approach is followed in~\cite{Spliet2015_discrete}, with the only difference that the first-stage time windows are selected from a finite set of \textit{a priori} constructed windows; this problem is referred to as the \textit{discrete} TWAVRP to distinguish it from the original \textit{continuous} TWAVRP.
In this paper, we consider both cases, and we shall in fact allow also for the generalized case in which feasible time window assignments lie in a continuous set for some portion of the customer base and in a discrete set for the remaining portion.

Algorithms to solve the aforementioned stochastic programming models have been proposed in~\cite{Spliet2015_continuous,Dalmeijer2018} for the continuous version, and in~\cite{Spliet2015_discrete} for the discrete version of the problem. The algorithms of~\cite{Spliet2015_continuous,Spliet2015_discrete} are based on branch-price-and-cut and can solve instances with 25 customers and 3 demand scenarios to optimality, while the algorithm of~\cite{Dalmeijer2018} is based on branch-and-cut and can address instances containing 40 customers and 3 scenarios. Several heuristics have also been proposed in~\cite{Spliet2015_discrete} for the discrete setting that can address instances containing up to 60 customers. Recently, \cite{Spliet2017} studied a variant of the TWAVRP with time-dependent travel times and proposed a branch-price-and-cut algorithm that can solve instances with 25 customers and 3 demand scenarios.

A problem that is closely related to the strategic TWAVRP is the Consistent Vehicle Routing Problem (ConVRP)~\cite{Groer2009}, which is  motivated in the context of operational level planning.
The ConVRP aims to design minimum cost vehicle routes over a finite, multi-day horizon to serve a set of customers with known demands.
The goal is to design routes that  are \emph{consistent} over time; this translates to satisfying any of the following requirements each time service is provided to a customer: \textit{(i)} arrival-time consistency, wherein the customer should be visited at roughly the same time during the day, \textit{(ii)} person-oriented consistency, in which the customer should be visited by the same driver, and whenever applicable, \textit{(iii)} load consistency, for which a customer should receive roughly the same quantity of goods. We refer the reader to~\cite{Kovacs2014} for an overview of this problem and its applications. 

Conceptually, the assigned time windows in the TWAVRP (which are also referred to as the \textit{endogenous time windows}) serve to satisfy the arrival-time consistency requirement of the ConVRP, which requires that every customer be visited at roughly the same time whenever service is requested.
Formally, the ConVRP requires that the difference between the earliest and the latest arrival times at each customer location must differ by no more than some pre-specified constant bound, which is referred to as the \textit{maximum allowable arrival-time differential}.
This bound is analogous to the pre-specified width of the endogenous time window in the continuous TWAVRP, and this equivalence between the two problems has been previously acknowledged in~\cite{Spliet2017,Spliet2015_continuous}.

The equivalence between the TWAVRP and the arrival-time ConVRP has two important consequences.
First, we observe that, in the most general case, the ConVRP allows for the possibility that not all customers require service in all time periods and that operational parameters (such as customer demands or travel times) differ from one time period to the other.
Translated in the context of the TWAVRP, this allows us to address applications in which a fraction of the customer base does not require frequent (e.g., daily) service (by considering scenarios where certain subsets of customers have no demand), as well as to treat uncertainty in a wider variety of parameters such as travel and service times (by considering scenarios in which their values represent perturbations from some nominal value). 
We therefore study the TWAVRP under a more general definition than what has been previously considered in the literature, in which it is possible to incorporate uncertainty in several operational parameters at once. However, we do remark that, as is the case with traditional stochastic programming models, the simultaneous treatment of uncertainty in several parameters may come at the cost of an explosion in the number of scenarios that have to be considered.

The second consequence of the equivalence between the TWAVRP and the arrival-time ConVRP is that any algorithm developed for the latter can be used to obtain solutions for the former.
In this paper, we adapt a decomposition algorithm proposed in~\cite{Subramanyam2017} for the Consistent Traveling Salesman Problem (ConTSP), the single-vehicle variant of the ConVRP that focuses purely on the aspect of arrival-time consistency, to obtain a new algorithm for the TWAVRP.
Our method can be viewed as a \textit{scenario decomposition algorithm} in the language of stochastic programming, and is not based on branch(-price)-and-cut that has been the de facto approach for solving TWAVRP models.
Our algorithm has the attractive features of \textit{modularity} and \textit{scalability}.
It is modular in accommodating \textit{(i)} continuous and discrete time windows, \textit{(ii)} any VRPTW solver (exact or heuristic), \textit{(iii)} routing-specific constraints (e.g., heterogeneous fleets), and \textit{(iv)} generic scenario descriptions.
Moreover, it can be readily parallelized which allows postulating a large number of scenarios of the uncertainty.
Together with \textit{(ii)} above, this means that our algorithm is also scalable.
The distinct contributions of our work may be summarized as follows.
\begin{itemize}
\item We generalize the definition of the TWAVRP. In particular, we study problems with scenario-based models of uncertainty in which any operational parameter may be uncertain and in which the endogenous time windows may be chosen from either continuous or discrete sets.

\item We propose an exact, scenario decomposition algorithm for the TWAVRP. The algorithm outperforms all state of the art methods, solving 54 out of 81 previously open instances. Furthermore, it can be easily parallelized, utilize any available VRPTW solver in a ``black box'' fashion, and be readily modified as a heuristic to solve large-scale instances.

\item We propose a new class of path-based disjunctions for the TWAVRP. We show via numerical experiments that these disjunctions, which manifest as branching rules in our algorithm, significantly improve the performance of the latter. In view of the relationship between the TWAVRP and ConVRP, these disjunctions are new for the ConVRP as well.

\item We conduct experiments with a parallel implementation of our algorithm to solve instances consisting of up to fifteen scenarios, representing a five-fold increase compared to existing literature. We use these solutions to elucidate, via out-of-sample simulations, the cost savings that are to be expected when considering more scenarios during time window assignment.
\end{itemize}

The rest of this paper is organized as follows. Section~\ref{sec:literature_review} reviews the relevant literature; in Section~\ref{sec:problem_definition}, we provide a general mathematical definition of the TWAVRP; in Section~\ref{sec:solution_approach}, we describe our solution algorithm for this problem; Section~\ref{sec:implementation_details} elaborates on important implementation details of our algorithm; Section~\ref{sec:computational_results} presents computational results on existing as well as new datasets; and, finally, we conclude in Section~\ref{sec:conclusions}.

\section{Related Literature}\label{sec:literature_review}
This section reviews papers in the vehicle routing literature, other than those mentioned in the introduction, which deal with aspects of \textit{(i)} endogenously imposed time windows, \textit{(ii)} consistent service considerations, and \textit{(iii)} stochastic or uncertain parameters. We choose not to review the extensive literature on the VRPTW; instead, we refer interested readers to~\cite{chapterInTothVigo}.

The authors of~\cite{Jabali2015} study a time window assignment problem that is encountered by courier companies who must quote delivery time windows to customers receiving sensitive packages. In this problem, which they refer to as the VRP with Self-Imposed Time Windows, travel times are uncertain but all customers and their demands are known \textit{a priori}. Using this information, the service provider must simultaneously determine \textit{(i)} a single routing plan to serve all customers, and \textit{(ii)} time window assignments that will be quoted to the customers before the vehicles depart from the depot. The objective is to minimize the sum of (deterministic) routing costs and (expected) overtime and tardiness penalty costs. Since travel times are uncertain, the key challenge is to determine the optimal \textit{placement} of time windows (along each vehicle route) in step \textit{(ii)} so as to avoid penalties due to delays. The uncertainty in travel time along each arc is modeled via a discrete set of ``disruption'' scenarios (each representing a deviation from some nominal value). Under the assumption that at most one arc will be disrupted on any vehicle route, the authors propose a Tabu Search heuristic for route generation and a linear programming approach for time window placement that inserts ``time buffers'' along each vehicle route. Recently, with a goal to solve the same problem, the authors of \cite{Vareias2017} extended the work of \cite{Jabali2015}. On the one hand, they relax the assumption that only one arc will be disrupted and use probabilistic chance constraints to guarantee reliable service. On the other hand, they propose an alternative model of uncertainty in which the stochastic deviations in travel times are modeled as continuous gamma-distributed random variables. Finally, by considering also the width of the time window (along with its placement) as a decision variable, they propose an Adaptive Large Neighborhood Search solution procedure.

A related line of work is the so-called Time Slot Management Problem that is motivated in the context of attended home delivery in e-commerce businesses~\cite{Agatz2011}. Here, customers place online orders for products (e.g., groceries) and, during this process, they select one time window (amongst a number of available ones) in which they want their product to be delivered. From the service provider's point of view, the challenge is to design a finite set of time windows (instead of just one) to offer to potential customers in different zip code areas. The problem is complicated by the fact that, during the design phase, the set of customers as well as their demand is not known with certainty. The objective is to design time windows that would not only yield low distribution costs in the long run, but also satisfy marketing or regulatory requirements. 
In general, existing approaches (e.g., see~\cite{Agatz2011,Hernandez2017,Bruck2017}) deal with uncertainty by simply using expected values of customer demand whose temporal distribution is either assumed to be uniform over the offered time windows or determined via simulations. The expected routing cost associated with a candidate set of time windows is then estimated via coarse continuous approximation methods (e.g., see~\cite{Daganzo1984,Figliozzi2009}) or detailed vehicle routing models. These estimates are embedded within some heuristic procedure (e.g., local search) to determine the final set of time windows.
We refer the reader to~\cite{Agatz2008} for an overview of research problems in the area of attended home delivery, including time slot management.
Finally, we mention the work of~\cite{Ulmer2017} who also study a time window assignment problem that is motivated in the context of home-attended services (e.g., cable installation). Here, as is the case in attended home delivery, customers dynamically place orders for some service. However, instead of the customer choosing a time window from a number of available ones, the service provider must quote a service time window to the customer at the time of request. Similar to the VRP with Self-Imposed Time Windows, travel and service times are stochastic and the objective is to minimize expected delays. However, unlike the latter problem, not all customers who will be serviced are known at the time when a particular request is received, and thus the customer base is also stochastic. The authors use approximate dynamic programming techniques to obtain time window assignments in real time.

We now review papers that study the ConVRP or its variants. The problem was introduced in \cite{Groer2009} motivated by operations encountered in small package shipping services offered by courier companies. Since then, there have been many studies that have focused on trying to design solution algorithms for the ConVRP (e.g., see~\cite{Groer2009,Tarantilis2012,KovacsALNS,Lian2017}). To the best of our knowledge, most of these approaches are heuristic in nature, and are based on the concept of generating a ``template'' routing plan that services only the most frequent customers; the daily routes are derived from the template by appropriately modifying it in a way that satisfies all service consistency requirements, in particular, driver and arrival-time consistency. With a focus on trying to address purely the requirement of arrival-time consistency, exact solution approaches for the ConTSP (the single-vehicle variant of the ConVRP) have been proposed in \cite{Subramanyam2016,Subramanyam2017}.
Given our result from later in this paper that establishes equivalence between the TWAVRP and arrival-time ConVRP, one can, in principle, adapt any of the aforementioned methods to obtain algorithms for the TWAVRP by ignoring aspects of driver consistency, whenever applicable. We, however, shall choose the decomposition algorithm of \cite{Subramanyam2017} for this purpose, due to a number of reasons. First, because of its decomposition principle, it solves the original problem by breaking it down into period-specific routing problems with time windows. This has the advantage that users can use their own routing solver as a module to obtain time window assignments and that it allows easy parallelization to solve instances with many time periods (or scenarios in the case of the TWAVRP). Second, it can be run in either exact or heuristic mode by changing the underlying routing solver. Third, as we will demonstrate, it can be readily extended to solve also the discrete TWAVRP, for which other ConVRP approaches would not be directly applicable. 

We conclude our literature review by mentioning the relationship of the TWAVRP to Stochastic Vehicle Routing Problems (SVRP). The latter class of problems also treats parameter uncertainty in the context of vehicle routing. However, unlike the TWAVRP which is inherently a strategic decision-making problem, the SVRP is an operational problem. Specifically, in the TWAVRP, the exact values of all parameters are assumed to be known before the vehicle routes are to be determined on a particular day. In contrast, in the SVRP, the vehicle routes must be determined before the parameter values become known, which are only gradually revealed during the execution of the routing plan. This requires fundamentally different modeling considerations and corresponding solution approaches.
The most common modeling paradigms are \textit{(i)} recourse models, \textit{(ii)} reoptimization models, and \textit{(iii)} chance-constrained models. 
In \textit{(i)}, a planned solution is designed in the first stage and recourse actions based on a predetermined policy are taken in the second stage when the uncertainties are revealed. For example, the capacity of a vehicle may get exceeded \textit{en route}, if demands are stochastic at the time of vehicle dispatch; in such cases, a recourse policy, such as a detour to the depot to empty the vehicle, must be explicitly incorporated in the model~\cite{Dror1989,Gauvin2014}.
In \textit{(ii)}, the planned solution is dynamically modified as the uncertain parameters (e.g., demands or travel times) become gradually revealed during the execution of the vehicle routes~\cite{Secomandi2001}.
Finally, in \textit{(iii)}, probabilistic or chance constraints are used to explicitly control the level of risk that is acceptable to the decision-maker~\cite{Laporte1989}.
We refer interested readers to~\cite{Gendreau2014}, who provide an excellent overview of applications, models and solution algorithms for the SVRP and its variants.

\section{Notation and Problem Definition}\label{sec:problem_definition}
\begin{figure}[!htb]
\centering
\begin{tikzpicture}[scale=0.7,%
					tw box/.style = {},%
					rarrow/.style = {-{>[width=2mm,length=2mm]}},%
					lrarrow/.style = {{<[width=2mm,length=2mm]}-{>[width=2mm,length=2mm]}}]
\begin{scope}
%\draw[help lines] (0, 0) grid (9, 4);
\draw [rarrow] (0, 0) -- +(9, 0) node[right] {\it time};
\draw (0.5, +.2) -- +(0, -.4) node[below] {$\vphantom{\ell_i}e_i$};
\draw (2.5, +.2) -- +(0, -.4) node[below] {$\vphantom{\ell_i}y_{i}$};
\draw (6.5, +.2) -- +(0, -.4) node[below] {$\vphantom{\ell_i}y_i + w_i$};
\draw (8.5, +.2) -- +(0, -.4) node[below] {$\ell_i$};

\draw [tw box] (2.5, 0.5) rectangle +(4, 1);

\draw [dashed] (2.5, 0) -- +(0, 2);
\draw [dashed] (6.5, 0) -- +(0, 2);

\draw [lrarrow] (2.5, 1.8) -- node[above] {$w_i$} +(4, 0);
\end{scope}

\begin{scope}[xshift = 12cm]
%\draw[help lines] (0, -.5) grid (9, 4.5);
\draw [rarrow] (0, 0) -- +(9, 0) node[right] {\it time};
\draw (0.5, +.2) -- +(0, -.4) node[below] {$\vphantom{\ell_i}e_i = \ubar{y}_{i1}$};
\draw (2.5, +.2) -- +(0, -.4) node[below] {$\vphantom{\ell_i}\bar{y}_{i1}$};
\draw (3.5, +.2) -- +(0, -.4) node[below] {$\vphantom{\ell_i}\ubar{y}_{i2}$};
\draw (5.5, +.2) -- +(0, -.4) node[below] {$\vphantom{\ell_i}\bar{y}_{i2}$};
\draw (6.5, +.2) -- +(0, -.4) node[below] {$\vphantom{\ell_i}\ubar{y}_{i3}$};
\draw (8.5, +.2) -- +(0, -.4) node[below] {$\bar{y}_{i3} = \ell_i$};

\draw [tw box] (0.5, 0.5) rectangle +(2, 1);
\draw [tw box] (3.5, 0.5) rectangle +(2, 1);
\draw [tw box] (6.5, 0.5) rectangle +(2, 1);

\draw [dashed] (0.5, 0) -- +(0, 2);
\draw [dashed] (2.5, 0) -- +(0, 2);
\draw [dashed] (3.5, 0) -- +(0, 2);
\draw [dashed] (5.5, 0) -- +(0, 2);
\draw [dashed] (6.5, 0) -- +(0, 2);
\draw [dashed] (8.5, 0) -- +(0, 2);

\end{scope}
\end{tikzpicture}
\caption{Illustration of continuous (left) and discrete (right) time window sets $TW_i$. In the continuous case, the assigned time window can be any sub-interval of $[e_i, \ell_i]$ of length $w_i$, while in the discrete case, the assigned window must be one of the intervals $[\ubar{y}_{i1}, \bar{y}_{i1}]$, $[\ubar{y}_{i2}, \bar{y}_{i2}]$ or $[\ubar{y}_{i3}, \bar{y}_{i3}]$.}
\label{figure:illustration_feasible_time_window_sets}
\end{figure}
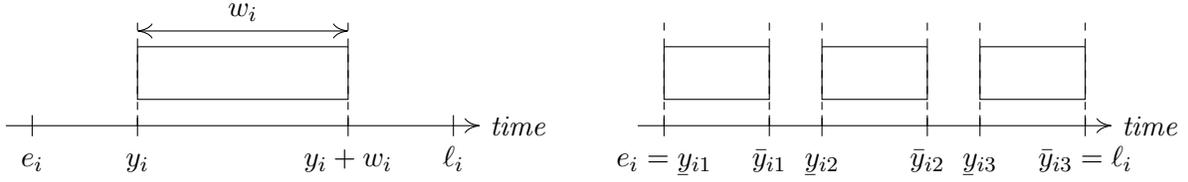

Let $G=(V,A)$ denote a directed graph with nodes $V = \{0,1,\ldots,n\}$ and arcs $A$. Node $0\in V$ represents the unique depot, and each node $i\in V_C \coloneqq V\setminus\{0\}$ represents a customer.
The operating hours of the depot are represented by the time window $[e_0, \ell_0]$, where an unlimited number of vehicles, each of capacity $Q$, are available for service. Each vehicle incurs a transportation cost $c_{ij}\in \mathbb{R}_{\geq 0}$ and a travel time $t_{ij}\in \mathbb{R}_{\geq 0}$ if it traverses the arc $(i, j) \in A$.
Furthermore, each customer $i \in V_C$ features a demand $q_i\in\mathbb{R}_{\geq 0}$, service time $u_i \in \mathbb{R}_{\geq 0}$ and \emph{exogenous} time window $[e_i, \ell_i]$ (e.g., representing operating hours).
The key decisions in the TWAVRP are to decide the \emph{endogenous} time windows $\tau_i \in TW_i$ to be assigned to each customer $i \in V_C$. The definition of the feasible time window set $TW_i$ may be either of the following (refer to Figure~\ref{figure:illustration_feasible_time_window_sets}):
\begin{itemize}
\item In the \emph{continuous} setting, the assigned time window must have a pre-specified width $w_i \in \mathbb{R}_{\geq 0}$; that is, $TW_i = \left\{[y_i, y_i + w_i] : e_i \leq y_i \leq \ell_i - w_i \right\}$, where we assume, without loss of generality, that $e_i \leq \ell_i - w_i$.

\item In the \emph{discrete} setting, the assigned time windows must belong to a pre-specified finite set; that is, $TW_i = \left\{[\ubar{y}_{i1}, \bar{y}_{i1}], \ldots, [\ubar{y}_{iN_i}, \bar{y}_{iN_i}]\right\}$, where we can assume, without loss of generality, that all $N_i$ candidate time windows are pairwise either non-overlapping or partially overlapping.\footnote{Two completely overlapping time windows $[a, b]$ and $[c, d]$ with $a \leq c \leq d \leq b$ can be replaced with the larger of the two time windows $[a, b]$.}
Therefore, the set $TW_i$ can be ordered so that $e_i = \ubar{y}_{i1} \leq \ldots \leq \ubar{y}_{iN_i}$ and $\bar{y}_{i1} \leq \ldots \leq \bar{y}_{iN_i} = \ell_i$.\footnote{We remark that $e_i = \ubar{y}_{i1}$ and $\ell_i = \bar{y}_{iN_i}$ can be achieved by preprocessing. If $e_i < \ubar{y}_{i1}$, then $e_i$ can be shifted forward to match $\ubar{y}_{i1}$, and if $e_i > \ubar{y}_{ib}$, then $\ubar{y}_{ib}$ can be shifted forward to match $e_i$. A similar argument applies for $\ell_i$.}
\end{itemize}
We denote by $V_\text{cont} \subseteq V_C$ and $V_\text{disc} \subseteq V_C$ the subset of customers whose feasible time window sets are continuous and discrete respectively. We note that $V_\text{cont} \cap V_\text{disc} = \emptyset$ and $V_\text{cont} \cup V_\text{disc} = V_C$.

In practice, operational parameters such as those related to the transportation network (costs $c$, travel times $t$) or the customers (demands $q$, service times $u$) are often not known with certainty at the strategic level when time windows must be allocated. 
Let $\theta$ denote the set of all operational parameters, and let $\mathbb{P}$ denote the joint probability distribution of $\theta$.
The goal of the TWAVRP is to assign the time windows $\tau_i \in TW_i$ in a way that minimizes the expected cost of routing:
\begin{equation}\label{eq:twavrp_stochastic_programming_formulation}
\begin{array}{r@{\;\;}l}
\displaystyle \mathop{\text{minimize}}_{\tau} & \displaystyle \mathbb{E}_{\theta \sim \mathbb{P}}\left[\mathtt{VRPTW}(\tau; \theta)\right] \\
\text{subject to} & \tau_i \in TW_i \quad \forall i \in V_C.
\end{array}
\end{equation}
In the above \emph{stochastic programming} formulation, $\mathtt{VRPTW}(\tau; \theta)$ denotes the minimum cost of the vehicle routing problem with time windows $\tau_i$, $i\in V_C$ and operational parameters $\theta$. 
A formal mathematical definition of $\mathtt{VRPTW}(\tau; \theta)$ follows.

\subsection{Mathematical Definition of $\mathtt{VRPTW}(\tau; \theta)$}\label{sec:problem_definition_VRPTW}
In this section, the time windows $\tau$ and operational parameters $\theta$ are assumed to be fixed to certain values in their domain and support respectively.
As far as the routing operation is concerned, we shall assume that each customer with non-zero demand must be visited exactly once by a single vehicle; that is, split deliveries are not allowed. 
In this regard, a \emph{route set} $\mathbf{R} = (R_1,\ldots,R_m)$, where $m \geq 1$, represents a partition of the customer set $V_C$. 
Here, $R_k = (R_{k,1},\ldots,R_{k,n_k})$ represents the $k^\text{th}$ vehicle route, $R_{k,l}$ the $l^\text{th}$ customer and $n_k$ the number of customers visited by vehicle $k$.
The cost of a route $R_k$ is evaluated as $c(R_k) = \sum_{l=0}^{n_k} c_{R_{k,l}, R_{k,l+1}}$, where we define $R_{k,0} = R_{k,n_k + 1} = 0$, and the cost of $\mathbf{R}$ is defined as $c(\mathbf{R}) = \sum_{k=1}^m c(R_k)$.
The route set $\mathbf{R}$ is feasible, if (i) all capacity constraints are satisfied, i.e., $\sum_{i\in R_k} q_i \leq Q$ for all $k \in \{1,\ldots,m\}$, and (ii) all time window constraints are satisfied, i.e., there exists a vector of arrival times, $\bm{a} \in \mathcal{X}(\mathbf{R}, \tau; \theta)$, where $\mathcal{X}(\mathbf{R}, \tau; \theta)$ is the feasible solution set of the following linear system of inequalities:
\begin{equation}\label{eq:arrival_time_definition}
\hspace{-0.1em}\mathcal{X}(\mathbf{R}, \tau; \theta)
=
\begin{aligned}\left\{
\bm{a} \in \mathbb{R}_{\geq 0}^n \left|
\begin{aligned}
& a_{R_{k,1}} \geq e_0 + t_{0,R_{k,1}}                                && \forall k \in K \coloneqq \{1,\ldots,m\} \\
& a_{R_{k,l+1}} - a_{R_{k,l}}\geq t_{R_{k,l},R_{k,l+1}} + u_{R_{k,l}} && \forall l\in \left\{1, \ldots, n_k - 1\right\}, \forall k \in K \\
& a_{R_{k,n_k}} \leq \ell_0 - t_{R_{k,n_k}, 0} - u_{R_{k,n_k}}        && \forall k \in K \\
& a_i \in \tau_i                                                      && \forall i \in V_C
\end{aligned}\right.\right\}
\end{aligned}
\end{equation}
In this definition, $a_{R_{k,l}}$ denotes the arrival time at location $R_{k,l}$, i.e., the arrival time at the $l^\text{th}$ location on the $k^\text{th}$ vehicle route.
The first three inequalities essentially require that the arrival time at any location must be at least as large as the sum of the arrival and service times in the previous location and the time to travel from the previous to the current location.
The last inequality requires the arrival time at customer location $i \in V_C$ to be within the time window $\tau_i$.
Observe that, by this definition, if a vehicle arrives at customer location $i\in V_C$ at a time earlier than $\ubar{\tau}_i \coloneqq \min_{t \in \tau_i} t$, then it is allowed to wait until $\ubar{\tau}_i$. However, arriving later than $\bar{\tau}_i \coloneqq \max_{t \in \tau_i} t $ is not permitted. We denote by $\mathcal{R}(\tau; \theta)$ the set of all feasible route sets for the given realization of operational parameters $\theta$ and time window assignment $\tau$. 
The value of $\mathtt{VRPTW}(\tau; \theta)$ can now be defined as the optimal value of the following optimization problem:
\begin{equation}\label{eq:model_vrptw}\tag{$\mathtt{VRPTW}(\tau; \theta)$}
\begin{array}{r@{\;\;}l}
\displaystyle \mathop{\text{minimize}}_{\mathbf{R}} & c\left(\mathbf{R}\right) \\
\text{subject to} & \mathbf{R}\in\mathcal{R}(\tau; \theta).
\end{array}
\end{equation}

\subsection{Deterministic Equivalent of Stochastic Programming Formulation}\label{sec:problem_definition_TWAVRP}
In practice, the exact joint probability distribution $\mathbb{P}$ is either not explicitly available or is hard to obtain. Indeed, even if it is known exactly, computing the objective function involves multi-dimensional integration of the function $\mathtt{VRPTW}(\tau; \cdot)$, which is practically impossible considering that the solution of the deterministic problem $\mathtt{VRPTW}(\tau; \theta)$ is itself challenging and cannot be obtained in closed form.
Instead, we assume that we are given a finite set of $S$ scenarios $\theta_1, \ldots, \theta_S$ along with associated probabilities of occurrence $p_1, \ldots, p_S$, where $p_s > 0$, $s \in \mathcal{S} \coloneqq \{1, \ldots, S\}$ and $\sum_{s \in \mathcal{S}} p_s = 1$.
In this situation, we seek to optimize the following \emph{deterministic equivalent} of problem~\eqref{eq:twavrp_stochastic_programming_formulation}, obtained by replacing the expectation with a sample average:
\begin{equation}\label{eq:twavrp_sample_average_formulation_conceptual}
\begin{array}{r@{\;\;}l}
\displaystyle \mathop{\text{minimize}}_{\tau} & \displaystyle \sum_{s \in \mathcal{S}} p_s \mathtt{VRPTW}(\tau; \theta_s) \\
\text{subject to} & \tau_i \in TW_i \quad \forall i \in V_C.
\end{array}
\end{equation}
Following the definition of~\ref{eq:model_vrptw}, the above sample average formulation~\eqref{eq:twavrp_sample_average_formulation_conceptual} can be equivalently represented as follows:
\begin{equation}\label{eq:twavrp_sample_average_formulation}
\begin{array}{r@{\;\;}l}
\displaystyle \mathop{\text{minimize}}_{\tau, \mathbf{R}} & \displaystyle \sum_{s \in \mathcal{S}} p_s c \left(\mathbf{R}_s \right) \\
\text{subject to} & \tau_i \in TW_i \quad \forall i \in V_C \\
& \mathbf{R}_s \in \mathcal{R}(\tau; \theta_s) \quad \forall s \in \mathcal{S}.
\end{array}
\end{equation}
The optimization problem~\eqref{eq:twavrp_sample_average_formulation} shall be our primary focus for the rest of the paper. We shall denote by $(\tau, \{{\mathbf{R}_s}\}_{s\in \mathcal{S}})$ a feasible solution to this problem.

\section{Solution Approach}\label{sec:solution_approach}
Our solution approach for the TWAVRP is motivated by the observation that, in the continuous setting (where $V_\text{disc} = \emptyset$), problem~\eqref{eq:twavrp_sample_average_formulation} can be reduced to an instance of the arrival-time ConVRP. 
Consequently, any algorithm to solve the latter class of problems can be used to solve problem~\eqref{eq:twavrp_sample_average_formulation}.
Section~\ref{sec:solution_approach_algorithm} presents an exact branch-and-bound algorithm for this purpose; we note, however, that the presented algorithm is more general-purpose, since it can also address the setting where $V_\text{disc} \neq \emptyset$.
Section~\ref{sec:solution_approach_new_disjunctions} presents new valid disjunctions that can be used as alternative branching rules in the algorithm; 
Section~\ref{sec:solution_approach_upper_bounds} presents upper bounding procedures (i.e., generating good time window assignments) in the context of our algorithm;
and finally, Section~\ref{sec:solution_approach_heuristic_modification} shows how the (exact) algorithm can be modified as a heuristic to solve large-scale instances.

\subsection{Overview of Exact Algorithm}\label{sec:solution_approach_algorithm}
We adapt the decomposition algorithm of~\cite{Subramanyam2017} developed for the Consistent Traveling Salesman Problem (the single-vehicle variant of the ConVRP) to solve the TWAVRP to optimality.
Notably, we extend the algorithm to incorporate also the case of discrete time windows.
This section presents the main ingredients of the algorithm translated into the TWAVRP context.

The algorithm uses a branch-and-bound tree search to identify the optimal time window assignments by solving within each node a set of VRPTW instances.
The tree is initialized with the original problem instance enforcing only the exogenous time windows $[e, \ell]$.
If valid time window assignments $\tau \in TW$ cannot be constructed using the optimal solution of the current node, the algorithm creates new nodes by using disjunctions~\eqref{eq:valid_disjunctions_continuous} and~\eqref{eq:valid_disjunctions_discrete} as branching rules.
The resulting branching rules are valid because, for every feasible solution $(\tau, \{\mathbf{R}_s\}_{s\in \mathcal{S}})$ in problem~\eqref{eq:twavrp_sample_average_formulation}, there exist arrival-time vectors $\bm{a}_s \in \mathcal{X}(\mathbf{R}_s, [e, \ell]; \theta_s)$ for each $s \in \mathcal{S}$ that satisfy disjunctions~\eqref{eq:valid_disjunctions}.
\begin{subequations}\label{eq:valid_disjunctions}
\begin{alignat}{2}
\left[a_{si} \leq \beta + w_i/2 \; \forall s\in\mathcal{S}\right] \;\vee \; & \left[a_{si} \geq \beta - w_i/2 \; \forall s\in\mathcal{S}\right] &\qquad \forall \beta\in \mathbb{R}, \; \forall i \in V_\text{cont} \label{eq:valid_disjunctions_continuous} \\
\left[a_{si} \leq \bar{y}_{ib} \; \forall s\in\mathcal{S}\right] \;\vee \; & \left[a_{si} \geq \ubar{y}_{i,b+1} \; \forall s\in\mathcal{S}\right] &\qquad \forall b \in \{1,\ldots,N_i - 1\}, \; \forall i \in V_\text{disc}. \label{eq:valid_disjunctions_discrete}
\end{alignat}
\end{subequations}

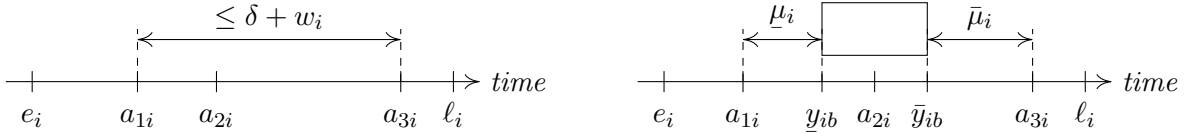
\begin{figure}[!htb]
\centering
\begin{tikzpicture}[scale=0.7,%
					tw box/.style = {},%
					rarrow/.style = {-{>[width=2mm,length=2mm]}},%
					lrarrow/.style = {{<[width=2mm,length=2mm]}-{>[width=2mm,length=2mm]}}]
\begin{scope}
%\draw[help lines] (0, 0) grid (9, 4);
\draw [rarrow] (0, 0) -- +(9, 0) node[right] {\it time};
\draw (0.5, +.2) -- +(0, -.4) node[below] {$\vphantom{\ell_i}e_i$};
\draw (2.5, +.2) -- +(0, -.4) node[below] {$\vphantom{\ell_i}a_{1i}$};
\draw (4.0, +.2) -- +(0, -.4) node[below] {$\vphantom{\ell_i}a_{2i}$};
\draw (7.5, +.2) -- +(0, -.4) node[below] {$\vphantom{\ell_i}a_{3i}$};
\draw (8.5, +.2) -- +(0, -.4) node[below] {$\ell_i$};

\draw [dashed] (2.5, 0) -- +(0, 1);
\draw [dashed] (7.5, 0) -- +(0, 1);

\draw [lrarrow] (2.5, 0.8) -- node[above] {$\leq \delta + w_i$} +(5, 0);
\end{scope}

\begin{scope}[xshift = 12cm]
%\draw[help lines] (0, -.5) grid (9, 4.5);
\draw [rarrow] (0, 0) -- +(9, 0) node[right] {\it time};
\draw (0.5, +.2) -- +(0, -.4) node[below] {$\vphantom{\ell_i}e_i$};
\draw (2.0, +.2) -- +(0, -.4) node[below] {$\vphantom{\ell_i}a_{1i}$};
\draw (3.5, +.2) -- +(0, -.4) node[below] {$\vphantom{\ell_i}\ubar{y}_{ib}$};
\draw (5.5, +.2) -- +(0, -.4) node[below] {$\vphantom{\ell_i}\bar{y}_{ib}$};
\draw (4.5, +.2) -- +(0, -.4) node[below] {$\vphantom{\ell_i}a_{2i}$};
\draw (7.5, +.2) -- +(0, -.4) node[below] {$\vphantom{\ell_i}a_{3i}$};
\draw (8.5, +.2) -- +(0, -.4) node[below] {$\ell_i$};

\draw [tw box] (3.5, 0.5) rectangle +(2, 1);

\draw [dashed] (2.0, 0) -- +(0, 1);
\draw [dashed] (3.5, 0) -- +(0, 1);
\draw [dashed] (5.5, 0) -- +(0, 1);
\draw [dashed] (7.5, 0) -- +(0, 1);

\draw [lrarrow] (2.0, 0.8) -- node[above] {$\ubar{\mu}_i$} +(1.5, 0);
\draw [lrarrow] (5.5, 0.8) -- node[above] {$\bar{\mu}_i$} +(2.0, 0);

\end{scope}
\end{tikzpicture}
\caption{Decision variables in problem~\eqref{eq:feasibility_milp} for the continuous (left) and discrete (right) cases.}
\label{figure:separation_problem}
\end{figure}

Conversely, if there exist route sets $\mathbf{R}_s \in \mathcal{R}([e, \ell]; \theta_s)$ and arrival-time vectors $\bm{a}_s \in \mathcal{X}(\mathbf{R}_s, [e, \ell]; \theta_s)$ for each $s \in \mathcal{S}$ satisfying disjunctions~\eqref{eq:valid_disjunctions}, then there exists a time window assignment $\tau \in TW$ such that $(\tau, \{\mathbf{R}_s\}_{s\in \mathcal{S}})$ is feasible in problem~\eqref{eq:twavrp_sample_average_formulation}. Specifically, a feasible time window assignment is
\begin{equation}\label{eq:actual_time_window_assignment}
\tau_i^\star = \begin{dcases}
\left[y_i, y_i + w_i\right], \text{where } y_i = \min\left\{ \ell_i - w_i, \min_{s\in\mathcal{S}} a_{si} \right\} & \text{if } i \in V_\text{cont} \\
\left[\ubar{y}_{ib_i}, \bar{y}_{ib_i}\right], \text{where } b_i \in \mathop{\arg\min}_{b \in \{1,\ldots,N_i\}} \left\{ \bar{y}_{ib}: \bar{y}_{ib} \geq \max_{s\in\mathcal{S}} a_{si} \right\} & \text{if } i \in V_\text{disc}
\end{dcases} \quad \forall i \in V_C.
\end{equation}
For given route sets $\mathbf{R}_s \in \mathcal{R}([e, \ell]; \theta_s)$, verifying the existence of arrival-time vectors $\bm{a}_s \in \mathcal{X}(\mathbf{R}_s, [e, \ell]; \theta_s)$, $s \in \mathcal{S}$, that satisfy disjunctions~\eqref{eq:valid_disjunctions} is equivalent to verifying that the optimal objective value $\delta$ of the following mixed-integer linear optimization problem is non-positive.

\begin{equation} \label{eq:feasibility_milp}
\begin{array}{l@{\quad}ll}
\displaystyle \mathop{\text{minimize}}_{\delta,a,z,\bar{\mu},\ubar{\mu}} & \delta \\
\text{subject to} & \displaystyle  \delta \in \mathbb{R}_{\geq 0}, \;\; \bm{a}_s \in \mathcal{X}(\mathbf{R}_s, \tau; \theta_s),\; s \in \mathcal{S} \\
& \left.\begin{array}{l}
\mspace{-11mu} \displaystyle \bar{\mu}_i, \ubar{\mu}_i \in \mathbb{R}_{\geq 0}, \;\; z_{ib} \in \{0, 1\},\; b \in \{1, \ldots, N_i\} \\
\mspace{-11mu} \displaystyle \sum_{b = 1}^{N_i} z_{ib} = 1 \\
\mspace{-11mu} \displaystyle z_{ib} = 1 \;\; \Rightarrow \;\; \bar{\mu}_i \geq a_{si} - \bar{y}_{ib}   \quad \forall (s, b) \in \mathcal{S} \times \{1, \ldots, N_i\} \\
\mspace{-11mu} \displaystyle z_{ib} = 1 \;\; \Rightarrow \;\; \ubar{\mu}_i \geq \ubar{y}_{ib} - a_{si} \quad \forall (s, b) \in \mathcal{S} \times \{1, \ldots, N_i\} 
\end{array}
\right\} & \forall i \in V_\text{disc} \\
& \displaystyle  \delta \geq a_{s_1i} - a_{s_2i} - w_i \quad \forall (s_1, s_2) \in \mathcal{S}\times\mathcal{S} & \forall i \in V_\text{cont} \\
& \displaystyle  \delta \geq \sum_{i \in V_\text{disc}} \left(\bar{\mu}_i + \ubar{\mu}_i \right).
\end{array}
\end{equation}
In this problem, $\delta$ records the minimum possible violation of disjunctions~\eqref{eq:valid_disjunctions} across all feasible arrival-time vectors $\bm{a}_s \in \mathcal{X}(\mathbf{R}_s, [e, \ell]; \theta_s)$, $s \in \mathcal{S}$.
The second-to-last constraint ensures that $\delta$ is at least as large as the maximum violation across all members of~\eqref{eq:valid_disjunctions_continuous}, while the last constraint ensures that $\delta$ is at least as large as the sum of violations across all members of~\eqref{eq:valid_disjunctions_discrete}.
In the latter case, the binary variable $z_{ib}$ indicates whether the $b^\text{th}$ member of~\eqref{eq:valid_disjunctions_discrete} is minimally violated for given $i \in V_\text{disc}$.
In other words, if $z_{ib} = 1$, then $\left[\ubar{y}_{ib}, \bar{y}_{ib}\right]$ is the best time window for customer~$i$, and $\bar{\mu}_i = \left[\max_{s \in \mathcal{S}} a_{si} - \bar{y}_{ib} \right]_{+}$ and $\ubar{\mu}_i = \left[\ubar{y}_{ib} - \min_{s \in \mathcal{S}} a_{si} \right]_{+}$ respectively record the arrival-time violations with respect to the start and end of this time window (refer to Figure~\ref{figure:separation_problem}). Here, $[ \cdot ]_+ = \max \{ \cdot, 0 \}$.

\paragraph{Algorithm.}
The above observations suggest that we can solve the TWAVRP without having to explicitly encode its time window assignments.
In fact, any TWAVRP instance decomposes into its individual scenarios if we track, within each node of a branch-and-bound search tree, a vector of applicable time windows (one for each customer), which we shall denote by $\tau$.
The root node enforces only the exogenous time windows (see Step~\ref{main_algo:init}).
Processing a node amounts to solving a set of VRPTW instances (one for each scenario) with time window constraints enforced by that node (see Step~\ref{main_algo:node_processing}).
It is important to remark that these VRPTW instances are uncoupled, and can thus be solved independently of each other.
This is because none of the expressions within each disjunct in~\eqref{eq:valid_disjunctions} (upon which our branching rules are based) link arrival-times from different scenarios in the same inequality.
The resulting optimal route sets $\{\mathbf{R}_s\}_{s \in \mathcal{S}}$ are then used as inputs to the separation problem~\eqref{eq:feasibility_milp} (see Step~\ref{main_algo:check_feasibility}).
If the optimal objective value satisfies $\delta \leq 0$, then a new, improved time window assignment $\tau^\star$ is recorded, as per~\eqref{eq:actual_time_window_assignment}.
Otherwise, the algorithm creates two new nodes with tightened time windows for customer $i^\star \in V_C$ (see Step~\ref{main_algo:branch}).
The overall algorithm is as follows.
\begin{enumerate}
\item\label{main_algo:init}
\textit{Initialize.} 
Set root node $\tau^0 \gets \left(\left[e_1, \ell_1\right], \ldots, \left[e_n, \ell_n\right]\right)$, node queue $\mathcal{N} \gets \left\{ \tau^0 \right\}$, upper bound $UB \gets +\infty$ and optimal time window assignment $\tau^\star \gets \emptyset$.

\item\label{main_algo:convergence}
\textit{Check convergence.} If $\mathcal{N} = \emptyset$, then stop: $\tau^\star$ is the optimal time window assignment with (expected) cost $UB$.
Otherwise, select a node $\tau$ from $\mathcal{N}$, and set $\mathcal{N} \gets \mathcal{N}\setminus \{\tau\}$.

\item\label{main_algo:node_processing}
\textit{Process node.} For each $s \in \mathcal{S}$, solve $\mathtt{VRPTW}(\tau; \theta_s)$; and let $\mathbf{R}_s$ denote an optimal solution.

\item\label{main_algo:fathom_by_bound}
\textit{Fathom by bound.} If $\sum_{s \in \mathcal{S}} p_s c(\mathbf{R}_s) \geq UB$, then go to Step~\ref{main_algo:convergence}.

\item\label{main_algo:check_feasibility}
\textit{Check feasibility.} Let $(\delta, a, z, \bar{\mu}, \ubar{\mu})$ be an optimal solution of problem~\eqref{eq:feasibility_milp}.
If $\delta \leq 0$, then a new, improved time window assignment is found: set $\tau^\star$ as per~\eqref{eq:actual_time_window_assignment}, set $UB \gets \sum_{s \in \mathcal{S}} p_s c(\mathbf{R}_s)$ and go to Step~\ref{main_algo:convergence}.

\item\label{main_algo:branch}
\textit{Branch.} Instantiate two children nodes, $\tau^\text{L}$ and $\tau^\text{R}$, from the parent node: $\tau^\text{L} \gets \tau$, $\tau^\text{R} \gets \tau$. If $\delta > \sum_{i \in V_\text{disc}} (\bar{\mu}_i + \ubar{\mu}_i)$, then do Step~\ref{main_algo:branch:continuous}; otherwise, do Step~\ref{main_algo:branch:discrete}:
\begin{enumerate}
\item\label{main_algo:branch:continuous}
\textit{Branch as per~\eqref{eq:valid_disjunctions_continuous}.}
Let $i^\star \in V_\text{cont}$ be any customer for which $\delta = \max_{s \in \mathcal{S}} a_{si^{\star}} - \min_{s \in \mathcal{S}} a_{si^\star} - w_{i^\star}$, and let $\beta^\star = \left(\max_{s \in \mathcal{S}} a_{si^{\star}} + \min_{s \in \mathcal{S}} a_{si^{\star}}\right)/2$. Tighten the time window for $i^\star$ as follows:
\textit{(i)} $\tau^\text{L}_{i^\star} \gets \left[\min_{t \in \tau_{i^\star}^\text{L}} t, \beta^\star + \frac12w_{i^\star} \right]$, 
\textit{(ii)} $\tau^\text{R}_{i^\star} \gets \left[\beta^\star - \frac12w_{i^\star}, \max_{t \in \tau_{i^\star}^\text{R}} t \right]$.  

\item\label{main_algo:branch:discrete}
\textit{Branch as per~\eqref{eq:valid_disjunctions_discrete}.}
Let $i^\star \in V_\text{disc}$ be any member of $\mathop{\arg\max}_{i \in V_\text{disc}} \{\bar{\mu}_i + \ubar{\mu}_i \}$. If $\bar{\mu}_{i^\star} \geq \ubar{\mu}_{i^\star}$, let $b^\star = \sum_{b = 1}^{N_i} b\mathbbm{1}_{\left[z_{ib} = 1\right]}$; else, let $b^\star = \sum_{b = 1}^{N_i} (b - 1)\mathbbm{1}_{\left[z_{ib} = 1\right]}$. 
Tighten the time window for $i^\star$ as follows:
\textit{(i)} $\tau^\text{L}_{i^\star} \gets \left[\min_{t \in \tau_{i^\star}^\text{L}} t, \bar{y}_{i^\star b^\star} \right]$, 
\textit{(ii)} $\tau^\text{R}_{i^\star} \gets \left[\ubar{y}_{i^\star\,b^\star + 1}, \max_{t \in \tau_{i^\star}^\text{R}} t \right]$.  
\end{enumerate}
Set $\mathcal{N}\gets \mathcal{N} \cup \{\tau^\text{L}, \tau^\text{R}\}$, and go to Step~\ref{main_algo:convergence}.
\end{enumerate}
We remark that any node selection strategy can be used in Step~\ref{main_algo:convergence} to guarantee convergence.

\begin{figure}[!htb]
\centering
\begin{tikzpicture}[scale=1.0,%
                    baseline={(current bounding box.center)},%
					tw box/.style = {},%
					edge/.style = {},%
					customer/.style = {circle, draw, minimum size = 0.5},%
					rarrow/.style = {-{>[width=2mm,length=2mm]}},%
					lrarrow/.style = {{<[width=2mm,length=2mm]}-{>[width=2mm,length=2mm]}}]
%\draw[help lines] (0, 0) grid (9, 4);
\node [customer] (3) at (0,0) {3};
\node [customer] (0) at (4,0) {0};
\node [customer] (1) at (6,0) {1};
\node [customer] (2) at (8,0) {2};
\node [customer] (4) at (5,2) {4};

\draw[edge] (3) -- node [below] {4} (0);
\draw[edge] (0) -- node [below] {2} (1);
\draw[edge] (1) -- node [below] {2} (2);
\draw[edge] (2) -- node [above] {6} (4);
\draw[edge] (1) -- node [right] {5} (4);
\draw[edge] (0) -- node [left]  {4} (4);
\draw[edge] (3) -- node [above] {7} (4);

\node[below = 0.1cm of 3] {$\left[0,15\right]$};
\node[below = 0.1cm of 0] {$\left[0,30\right]$};
\node[below = 0.1cm of 1] {$\left[0,10\right]$};
\node[below = 0.1cm of 2] {$\left[0,20\right]$};
\node[above = 0.1cm of 4] {$\left[0, 7\right], \left[14, 22\right]$};
\end{tikzpicture}$\quad$
\begin{tabular}{lccccc}
    \multicolumn{6}{l}{$(n, S, Q) = (4, 2, 3)$}\\
	\toprule
                 & $p_s$ & $q_{s1}$ & $q_{s2}$ & $q_{s3}$ & $q_{s4}$ \\ \midrule
    $s = 1\quad$ & 1/2   & 1        & 3        & 1        & 1 \\
    $s = 2\quad$ & 1/2   & 3        & 1        & 1        & 1 \\ \bottomrule
    \multicolumn{6}{l}{$w_1 = w_2 = w_3 = 3$}\\
\end{tabular}
\caption{Instance parameters for the illustrative example.
The arc weights denote both travel times and costs, and they are such that the triangle inequality is satisfied. Note that the graph is assumed to be complete, with nodes 0--3 lying along a straight line and costs/times to travel along this line being cumulative (not all arcs are shown for convenience). All service times are zero. The depicted intervals denote exogenous time windows for nodes 1--3 and feasible ones for node 4.}
\label{figure:example_instance}
\end{figure}
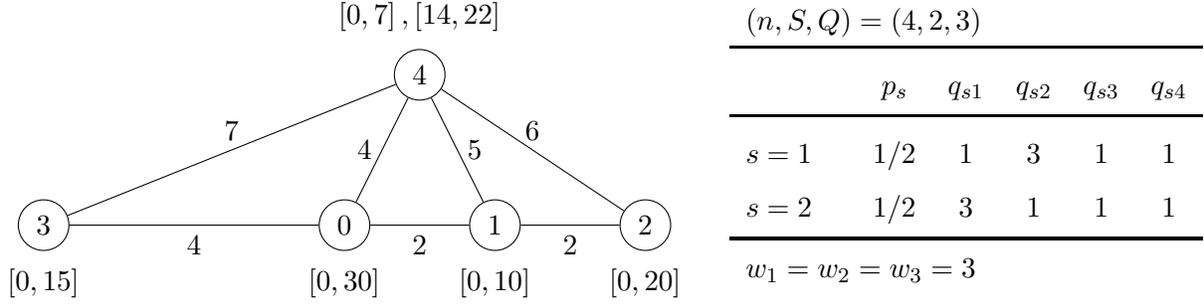

An illustration on a small example is now presented to aid understanding and give intuition about the algorithm. Consider the TWAVRP instance shown in Figure~\ref{figure:example_instance}. This example features $n= 4$ customers, with $V_\text{cont} = \{1,2,3\}$ and $V_\text{disc} = \{4\}$. Only customer demands are uncertain and they are represented using $S = 2$ scenarios.

\begin{figure}[!htb]
\centering
\input{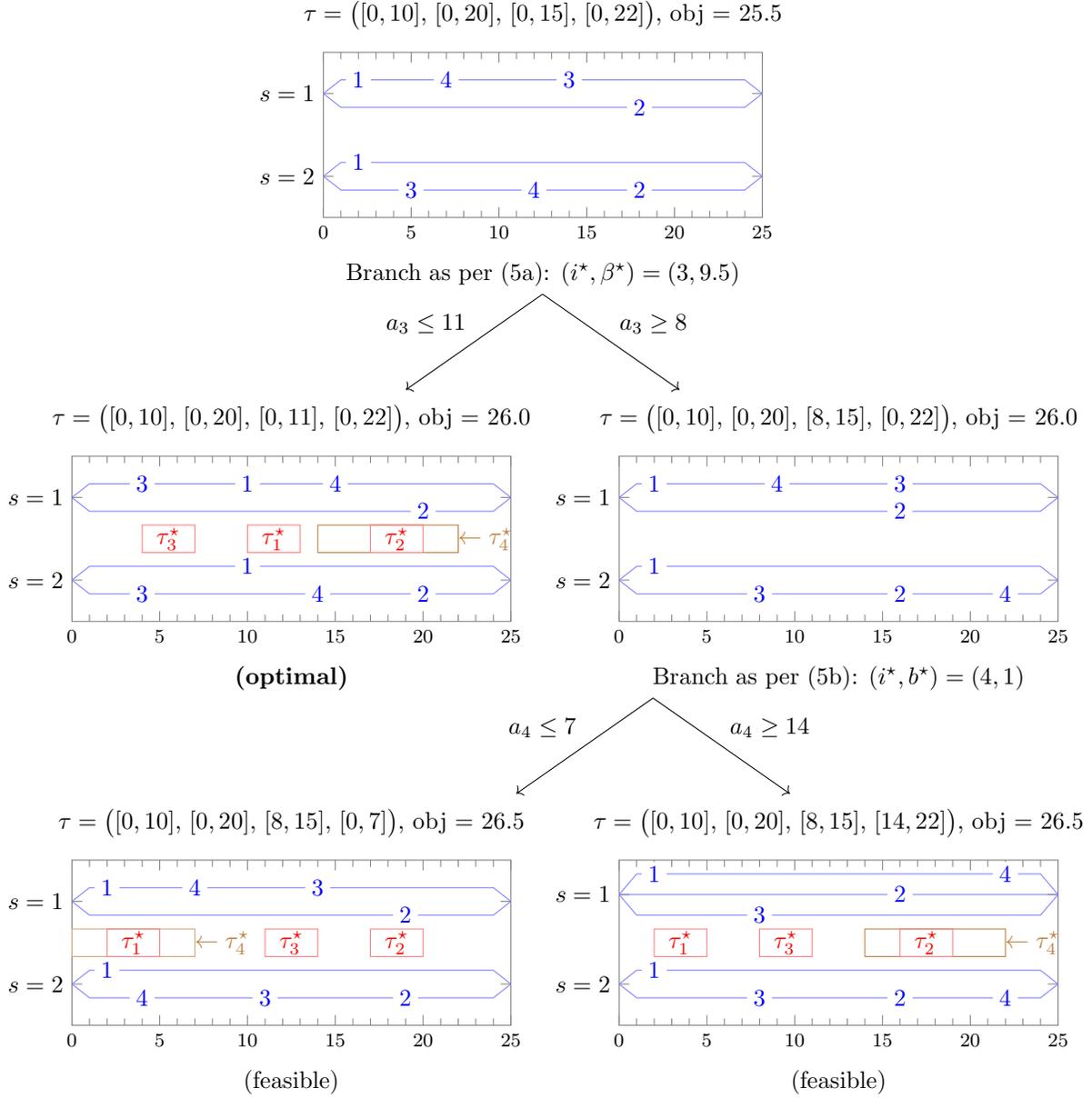}
\caption{The search tree of our algorithm for the illustrative example of Figure~\ref{figure:example_instance}.}
\label{figure:example_algorithm}
\end{figure}

The search tree of our algorithm to solve the illustrative example of Figure~\ref{figure:example_instance} is shown in Figure~\ref{figure:example_algorithm}. Each ``rectangle'' denotes a node of our search tree. Within each rectangle, for each scenario $s \in \{1, 2\}$, the optimal route set $\mathbf{R}_s$ is shown. The $x$-coordinate of each customer $i \in V_C$ denotes its arrival-time $a_{si}$, which is computed by solving the separation problem~\eqref{eq:feasibility_milp} to optimality (note that there may be multiple optimal solutions in each case). Finally, on top of each rectangle, the time windows $\tau$ enforced by our algorithm and the objective value of the corresponding optimal route sets (equal to $\sum_{s \in \mathcal{S}} p_s c(\mathbf{R}_s)$) are shown.
In the root node, the separation problem~\eqref{eq:feasibility_milp} certifies (in~Step~\ref{main_algo:check_feasibility}) that no valid time window assignment can be constructed from its optimal solution.
In particular, if we focus on customer $3 \in V_\text{cont}$, then $a_{13} = 14$ and $a_{23} = 5$. These arrival times clearly do not fall within a time window of width $w_3 = 3$.
Therefore, as per Step~\ref{main_algo:branch:continuous}, $\beta^\star = (14 + 5)/2 = 9.5$. The time window of customer $3$ in the left child is tightened to $[e_3, \beta^\star + w_3/2] = [0, 11]$, while in the right child to $[\beta^\star - w_3/2, \ell_3] = [8, 15]$.
Similarly, if we focus on customer $4 \in V_\text{disc}$ in this right child node, then $a_{14} = 9$ and $a_{24} = 22$. These arrival times also do not simultaneously satisfy either of the two candidate windows, $[0, 7]$ or $[14, 22]$. Therefore, a branch is made, as per Step~\ref{main_algo:branch:discrete}, tightening the time window of customer $4$ in the left child to $[0, 7]$ and in the right one to $[14, 22]$.

We remark that, in a given node of our search tree (except the root node), one does not need to solve a VRPTW subproblem for every scenario as required by Step~\ref{main_algo:node_processing} of the algorithm, and the optimal VRPTW route sets for some of the scenarios can be directly transferred from the parent subproblems.
In fact, after branching has occurred in Step~\ref{main_algo:branch} of the algorithm, at most $S$ (out of a total of $2S$) VRPTW subproblems have to be solved across the two children nodes.
This is because, by construction, any arrival time vectors corresponding to the optimal solution, $\mathbf{R}_s$, of a given scenario $s$ cannot simultaneously violate both disjuncts of the applied branching disjunction (either~\eqref{eq:valid_disjunctions_continuous} or~\eqref{eq:valid_disjunctions_discrete}), although it may satisfy both disjuncts simultaneously.
Therefore, as far as a given scenario $s$ is concerned, a VRPTW subproblem needs to be solved only at most once across the two children nodes.
To illustrate this, see Figure~\ref{figure:example_algorithm}. After branching has occurred in the root node, $\mathbf{R}_1$ in the right child node is exactly the same as that in the parent, since it already satisfies the applied branching constraint $[a_3 \geq 8]$. For the same reason, $\mathbf{R}_2$ in the left child node is exactly the same as that in the parent.
Once the branching rule has been established, whether a scenario-specific set of routes remains feasible (and hence, optimal) for the VRPTW instance of a child node can be inferred trivially by inspection, and hence, the corresponding instance need not be solved, as the routes can be copied over.
For the VRPTW instances that indeed warrant a new route set to be computed, Section~\ref{sec:implementation_details} describes methods for solving the corresponding subproblems.

\subsection{Path-based Disjunctions}\label{sec:solution_approach_new_disjunctions}
The algorithm described in the previous section converges in finite time (the argument is similar to that in~\cite{Subramanyam2017}).
However, the branching Step~\ref{main_algo:branch} may not necessarily ``cut off'' the VRPTW route set $\{ \mathbf{R}_s \}_{s \in \mathcal{S}}$ found in the parent node.
Indeed, it is only guaranteed to cut off the arrival time vector $\{\bm{a}_s\}_{s\in \mathcal{S}}$  corresponding to $\{ \mathbf{R}_s \}_{s \in \mathcal{S}}$.
More specifically, an optimal route set of the (right) child node $\tau^\text{R}$ may also be optimal for its parent node $\tau$. This is because the only difference between $\tau^\text{R}$ and its parent $\tau$ is that the former features a tighter earliest start time for customer $i^\star \in V_C$.
It is possible for an optimal route set of $\tau$ to be unchanged, if the tighter earliest start time constraint can be satisfied by simply \textit{allowing the vehicle to wait longer at $i^\star$}.
To see this, consider again the example of Figure~\ref{figure:example_instance}. If the exogenous time window of customer $2$ is increased by one unit to $[0, 21]$, then an optimal route set $\mathbf{R}_2$ in the root node of the algorithm (see Figure~\ref{figure:example_algorithm}) will also be optimal in its right child node since the vehicle visiting customer $3$ will simply wait longer until it satisfies the branching constraint, $a_3 \geq 8$. Along with the fact that $\mathbf{R}_1$ is also the same (refer to the discussion in the previous section), this means that the computed VRPTW route set $\{ \mathbf{R}_s \}_{s \in \mathcal{S}}$ in the root node has not changed in its right child node.
The impact of this is non-improving lower bounds: the objective value of the node $\tau^\text{R}$ will be exactly the same as that of its parent, leading to slow convergence and poor numerical performance.
This observation motivates us to investigate branching rules which are guaranteed to cut off the parent route set.

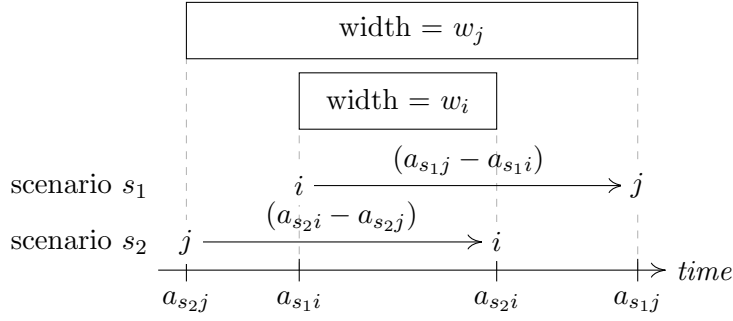
\begin{figure}[!htb]
\centering
\begin{tikzpicture}[scale=0.75,%
					tw box/.style = {},%
					projection/.style = {dashed,opacity=0.3},
					rarrow/.style = {-{>[width=2mm,length=2mm]}},%
					lrarrow/.style = {{<[width=2mm,length=2mm]}-{>[width=2mm,length=2mm]}}]
%\draw[help lines] (0, 0) grid (9, 4);
\draw [rarrow] (0, 0) -- +(9, 0) node[right] {\it time};
\draw (0.5, +.2) -- +(0, -.4) node[below] {$a_{s_2j}$};
\draw (2.5, +.2) -- +(0, -.4) node[below] {$a_{s_1i}$};
\draw (6.0, +.2) -- +(0, -.4) node[below] {$a_{s_2i}$};
\draw (8.5, +.2) -- +(0, -.4) node[below] {$a_{s_1j}$};

\draw [tw box] (2.5, 2.5) rectangle node {width = $w_i$} +(3.5, 1);
\draw [tw box] (0.5, 3.75) rectangle node {width = $w_j$} +(8, 1);

\draw [projection] (0.5, 0) -- +(0, 3.75);
\draw [projection] (2.5, 0) -- +(0, 2.5);
\draw [projection] (6.0, 0) -- +(0, 2.5);
\draw [projection] (8.5, 0) -- +(0, 3.75);

\draw [rarrow] (0.75, 0.5) -- node[above] {$(a_{s_2i} - a_{s_2j})$} +(5, 0);
\draw [rarrow] (2.75, 1.5) -- node[above] {$(a_{s_1j} - a_{s_1i})$} +(5.5, 0);
\node[fill=white,minimum size=0.5] at (0.5, 0.5) {$j$};
\node[fill=white,minimum size=0.5] at (6.0, 0.5) {$i$};
\node[fill=white,minimum size=0.5] at (2.5, 1.5) {$i$};
\node[fill=white,minimum size=0.5] at (8.5, 1.5) {$j$};
\node[left] at (0, 1.5) {scenario $s_1$};
\node[left] at (0, 0.5) {scenario $s_2$};
\end{tikzpicture}
\caption{Motivation for the path-based disjunctions. The above shown sub-paths must satisfy $(a_{s_2i} - a_{s_2j}) + (a_{s_1j} - a_{s_1i}) \leq w_i + w_j$, if they are to be part of a feasible solution.}
\label{figure:path_based_disjunction_motivation}
\end{figure}

Our motivation for the new class of disjunctions comes from the \emph{path precedence inequalities} proposed in~\cite{Dalmeijer2018} for the continuous TWAVRP and the \emph{inconsistent path elimination constraints} proposed in~\cite{Subramanyam2016} for the ConTSP.
Consider a feasible solution $(\tau, \{\mathbf{R}_s\}_{s \in \mathcal{S}})$ to problem~\eqref{eq:twavrp_sample_average_formulation}.
Suppose that there is a vehicle route in the solution $\mathbf{R}_{s_1}$ of scenario $s_1 \in \mathcal{S}$ in which customer $i \in V_C$ is visited before customer $j \in V_C \setminus \{i\}$, and that there is a vehicle route in the solution $\mathbf{R}_{s_2}$ of scenario $s_2 \in \mathcal{S} \setminus \{s_1\}$ in which $j$ is visited before $i$.
Since both $i$ and $j$ are visited within their respective time windows $\tau_i$ and $\tau_j$ in both scenarios, it must be the case that the sum of the travel times from $i$ to $j$ in scenario $s_1$ and from $j$ to $i$ in scenario $s_2$ is at most the sum of the widths of their time windows, $w_i + w_j$, as shown in Figure~\ref{figure:path_based_disjunction_motivation}.
Consequently, if this condition is not satisfied by a route set $\{\mathbf{R}_s\}_{s \in \mathcal{S}}$ for all possible pairs $(i,j)$, then there cannot exist a feasible time window assignment $\tau \in TW$ such that $(\tau, \{\mathbf{R}_s\}_{s \in \mathcal{S}})$ is feasible in problem~\eqref{eq:twavrp_sample_average_formulation}.

\begin{figure}[!htb]
\centering
\begin{tikzpicture}[scale=1,%
					tw/.style = {draw=red!50,text=red},%
					route/.style = {blue!50},%
					projection/.style = {draw=none,ultra thin},%
					customer/.style = {text=blue,circle, inner sep = 0pt, minimum size=0.35cm, fill=white}]
\small
\begin{axis}[xticklabels={,0,2,4,6,8},
			 xmin=0,xmax=8,
			 minor x tick num=1]
	\node[customer] (1s1) at (axis cs:2,10){1}; \draw[projection] (1s1) -- (axis cs:2,0);
	\node[customer] (2s1) at (axis cs:3,10){2}; \draw[projection] (2s1) -- (axis cs:3,0);
	\node[customer] (3s1) at (axis cs:4,10){3}; \draw[projection] (3s1) -- (axis cs:4,0);
	\node[customer] (4s1) at (axis cs:4,8) {4}; %\draw[projection] (4s1) -- (axis cs:7,0);
	\coordinate (route1s1start) at (axis cs:1,10);
	\coordinate (route2s1start) at (axis cs:1,8);
	\coordinate (route1s1end) at (axis cs:7,10);
	\coordinate (route2s1end) at (axis cs:7,8);
	\draw[route] (axis cs:0,9) -- (route1s1start) -- (1s1) -- (2s1) -- (3s1) -- (route1s1end) -- (axis cs:8,9);
	\draw[route] (axis cs:0,9) -- (route2s1start) -- (4s1) -- (route2s1end) -- (axis cs:8,9);
	
	\node[customer] (1s2) at (axis cs:6,4) {1}; \draw[projection] (1s2) -- (axis cs:6,0);
	\node[customer] (2s2) at (axis cs:5,4) {2}; \draw[projection] (2s2) -- (axis cs:5,0);
	\node[customer] (3s2) at (axis cs:4,2) {3}; %\draw[projection] (3s2) -- (axis cs:4,0);
	\node[customer] (4s2) at (axis cs:4,4) {4}; %\draw[projection] (4s2) -- (axis cs:4,0);
	\coordinate (route1s2start) at (axis cs:1,4);
	\coordinate (route2s2start) at (axis cs:1,2);
	\coordinate (route1s2end) at (axis cs:7,4);
	\coordinate (route2s2end) at (axis cs:7,2);
	\draw[route] (axis cs:0,3) -- (route1s2start) -- (4s2) -- (2s2) -- (1s2) -- (route1s2end) -- (axis cs:8,3);
	\draw[route] (axis cs:0,3) -- (route2s2start) -- (3s2) -- (route2s2end) -- (axis cs:8,3);
\end{axis}
\end{tikzpicture}
\caption{Counter-example to show that the disjunctions~\eqref{eq:path_disjunctions} are not sufficient for the TWAVRP.}
\label{figure:counterexample_path_disjunction}
\end{figure}
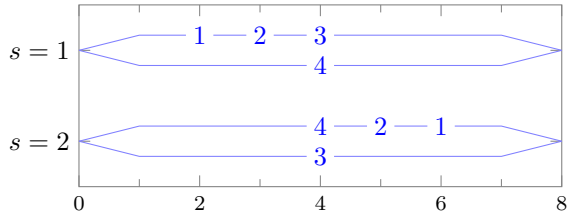

To construct valid disjunctions based on the above observation, we first introduce some notation.
Let $\pi = (v_1, \ldots, v_p)$ denote a directed $v_1 - v_p$ path in graph $G$ that is formed by the arcs in the set $\{(v_i, v_{i+1}): i = 1, \ldots, p-1 \}$, where $(v_i, v_{i+1}) \in A$ for all $i = 1,\ldots, p-1$.
We shall only consider paths which are open and simple, i.e., $p > 1$ and $v_i \neq v_j$ for $i \neq j$.
For a given realization of the travel and service times, the travel time along $\pi$ is defined to be $t(\pi) = \sum_{i = 1}^{p-1} \left( t_{v_{i}v_{i+1}} + u_{v_i} \right)$, where we define $u_0 = 0$. Note that as per this definition, the travel time along a path does not include any waiting time that might potentially be incurred at its nodes.
Finally, a route set $\mathbf{R} = (R_1, \ldots, R_m)$, where $R_k = (R_{k,1}, \ldots, R_{k,n_k})$ for $k = 1, \ldots, m$, is said to contain $\pi$, if $\pi$ appears as a sub-path in any of its routes, i.e., 
if for some $k$, we have $R_k = (R_{k,1}, \ldots, R_{k,l} = v_1, R_{k,l+1} = v_2, \ldots, R_{k,l+p-1} = v_p, \ldots, R_{k,n_k})$.
The key result of this section is that every feasible solution $(\tau, \{{\mathbf{R}_s}\}_{s\in \mathcal{S}})$ to problem~\eqref{eq:twavrp_sample_average_formulation} satisfies the following disjunctions:
\begin{equation}\label{eq:path_disjunctions}
\begin{aligned}
\left[\begin{array}{c}
\mathbf{R}_s \text{ does not contain any } i-j\\
\text{path with travel time} \geq d_1
\end{array} \;\; \forall s \in \mathcal{S} \right]
\vee
\left[\begin{array}{c}
\mathbf{R}_s \text{ does not contain any } j-i \\
\text{path with travel time} \geq d_2
\end{array} \;\; \forall s \in \mathcal{S} \right] \\
\forall (d_1, d_2) \in \mathbb{R}^2: d_1 + d_2 > w_i + w_j, \;\; \forall (i,j) \in V_C \times V_C: i\neq j,
\end{aligned}
\end{equation}
where $w_k$ for any $k \in V_\text{disc}$ is defined to be $w_k = \max_{b \in \{1, \ldots, N_k\}} \{\bar{y}_{kb} - \ubar{y}_{kb} \}$.
However, the converse is not true. 
To see this, consider the following counter-example.
\begin{itemize}[itemsep=0pt]
\item $(n, S, Q) = (4, 2, 3)$. $V_\text{cont} = V_C$ and $w_i = 1$ for all $i \in V_C$. 
\item $[e, \ell] = \left([0, 6], [0, 6], [3, 4], [4, 5]\right)$. Also, $[e_0, \ell_0] = [0, 10]$.
\item $G = (V, A)$ is complete. $c_{ij} = t_{ij} = 1$ for all $(i, j) \in A$ and $u_i = 0$ for all $i \in V_C$.
\item Demand is uncertain. $(q_{s1}, q_{s2}, q_{s3}, q_{s4}) = (1, 1, 1, 3)$ for $s=1$ and $(1, 1, 3, 1)$ for $s = 2$.
\end{itemize}
Consider the route sets $\{\mathbf{R}_s\}_{s=1,2}$ shown in Figure~\ref{figure:counterexample_path_disjunction}.
The $x$-coordinates correspond to arrival times.
This solution satisfies all path-based disjunctions~\eqref{eq:path_disjunctions}.
However, it is not a feasible TWAVRP solution since there are no valid time window assignments $\tau \in TW$ for customers~$1$ and~$2$.
In contrast, observe that this solution does indeed violate the (necessary and sufficient) time window-based disjunctions~\eqref{eq:valid_disjunctions_continuous} corresponding to $(i, \beta) = (1, 4)$ as well as $(i, \beta) = (2, 4)$.

The above observations suggest that we can use the path-based disjunctions~\eqref{eq:path_disjunctions} as the basis of a branching rule in addition to the time window-based disjunctions~\eqref{eq:valid_disjunctions}.
The corresponding branching constraints, i.e., constraints within each individual disjunct of~\eqref{eq:path_disjunctions}, are equivalent to \emph{path elimination constraints} (e.g., see~\cite{Ascheuer2001}).
Consequently, the subproblems to be solved in Step~\ref{main_algo:node_processing} of the algorithm are VRPTW instances with additional path elimination constraints.
To describe these subproblems formally, if $\mathcal{F}$ is a finite collection of triples $(i,j,d) \in V_C \times V_C \times \mathbb{R}$, such that $i\neq j$, then we let each member of $\mathcal{F}$ represent a family of \emph{forbidden paths}.
Specifically, the member $(i,j,d) \in \mathcal{F}$ represents the set of all $i-j$ paths with travel time greater than or equal to $d$. 
Given a collection $\mathcal{F}$ of representative forbidden paths, time windows $\tau$ and operational parameters $\theta$, we define $\mathtt{VRPTWFP}(\tau, \mathcal{F}; \theta)$ to be the optimal value of the following optimization problem:
\begin{equation}\label{eq:model_vrptw_with_path_elimination}\tag{$\mathtt{VRPTWFP}(\tau, \mathcal{F}; \theta)$}
\begin{array}{r@{\;\;}l}
\displaystyle \mathop{\text{minimize}}_{\mathbf{R}} & c\left(\mathbf{R}\right) \\
\text{subject to} & \mathbf{R}\in\mathcal{R}(\tau; \theta) \\
& \mathbf{R} \text{ does not contain any $i$-$j$ path } \pi \text{ with } t(\pi) \geq d  \quad \forall (i,j,d) \in \mathcal{F}.
\end{array}
\end{equation}
Note that $\mathtt{VRPTWFP}(\tau, \emptyset; \theta) = \mathtt{VRPTW}(\tau; \theta)$ so our definition is consistent with the one in Section~\ref{sec:problem_definition_VRPTW}.

Before we can incorporate the path-based disjunctions~\eqref{eq:path_disjunctions} as branching rules in our algorithm, we also need a separation algorithm, which will take as input route sets $\{\mathbf{R}_s\}_{s \in \mathcal{S}}$, and return either a violated member of~\eqref{eq:path_disjunctions} or a certificate that all of its members are satisfied. We compute the following quantities, where we assume $(i, j) \in V_C \times V_C$, such that $i \neq j$, is a given pair of customers.
\begin{itemize}[itemsep=0pt]
\item $\mathcal{S}_{ij}$: scenarios containing an $i-j$ path; that is, $\mathcal{S}_{ij} = \left\{ s \in \mathcal{S}: \mathbf{R}_s  \text{ contains an } i-j \text{ path} \right\}$.
\item $d_{ij}^s$: travel time of $i-j$ path in $\mathbf{R}_s$, where $s \in \mathcal{S}_{ij}$.
\item $\nu_{ij}$: no. of violating scenario pairs; that is, $\nu_{ij} = \left\lvert \left\{(s_1, s_2) \in \mathcal{S}_{ij} \times \mathcal{S}_{ji} : d_{ij}^{s_1} + d_{ji}^{s_2} > w_i + w_j \right\}\right\rvert$.
\item $\Delta_{ij}$: minimum value of sum of path travel times (across violating scenario pairs); that is, $\Delta_{ij} = \inf_{(s_1, s_2) \in \mathcal{S}_{ij} \times \mathcal{S}_{ji}} \left\{d_{ij}^{s_1} + d_{ji}^{s_2} : d_{ij}^{s_1} + d_{ji}^{s_2} > w_i + w_j\right\}$.
\end{itemize}
We are now in a position to incorporate the path-based-disjunctions~\eqref{eq:path_disjunctions} in our algorithm. To do so, we store the set of forbidden paths $\mathcal{F}$ as part of a node's characteristic data (along with $\tau$) at the time of node creation (initialization and branching steps), and we use this set as input to the VRPTW subproblem at the time of node processing.
We update $\mathcal{F}$ in the branching step based on a ranked list $\mathcal{L}$ of pairs $(i, j)$ for which the corresponding path-based disjunctions~\eqref{eq:path_disjunctions} can be used as branching rules.
In all, the following modifications are made to the main algorithm.
\begin{enumerate}
\item[\mylabel{path_based_algo:init}{\ref{main_algo:init}$^\prime$}] 
\textit{Initialize.} 
Set root node $\tau^0 \gets \left(\left[e_1, \ell_1\right], \ldots, \left[e_n, \ell_n\right]\right)$, $\mathcal{F}^0 \gets \emptyset$, node queue $\mathcal{N} \gets \left\{ (\tau^0, \mathcal{F}^0) \right\}$, upper bound $UB \gets +\infty$ and optimal time window assignment $\tau^\star \gets \emptyset$.

\item[\mylabel{path_based_algo:node_processing}{\ref{main_algo:node_processing}$^\prime$}] 
\textit{Process node.} For each $s \in \mathcal{S}$, solve $\mathtt{VRPTWFP}(\tau, \mathcal{F}; \theta_s)$; and let $\mathbf{R}_s$ denote its optimal solution.

\item[\mylabel{path_based_algo:branch}{\ref{main_algo:branch}$^\prime$}] 
\textit{Branch.} Instantiate two children nodes from $(\tau, \mathcal{F})$: $(\tau^\text{L}, \mathcal{F}^\text{L}) \gets (\tau, \mathcal{F})$,  $(\tau^\text{R}, \mathcal{F}^\text{R}) \gets (\tau, \mathcal{F})$.
Set $\mathcal{L} \gets \emptyset$. For each pair $(i, j) \in V_C \times V_C$, such that $i \neq j$ and $\nu_{ij} \geq 1$, set $\mathcal{L} \gets \mathcal{L} \cup \{(i, j)\}$.\\
If $\mathcal{L} \neq \emptyset$, do Step~6$^\prime$c; else if $\delta > \sum_{i \in V_\text{disc}} (\mu^{+}_i + \mu^{-}_i)$, do Step~\ref{main_algo:branch:continuous}; else, do Step~\ref{main_algo:branch:discrete}:
\begin{enumerate}
\item[(c)]\label{path_based_algo:branch:paths}
Sort $\mathcal{L}$ in decreasing order of $\nu$ (breaking ties in increasing order of $\Delta$).
Let $(i, j)$ be the first element of $\mathcal{L}$ and let $(s_1, s_2) = \mathop{\arg\min}_{(s, s^\prime) \in \mathcal{S}_{ij} \times \mathcal{S}_{ji}} \Big\{d_{ij}^{s} + d_{ji}^{s^\prime} : d_{ij}^{s} + d_{ji}^{s^\prime} > w_{i} + w_{j}\Big\}$.
Set $d_1^\star$ and $d_2^\star$ as follows: if $d_{ij}^{s_1} \leq d_{ji}^{s_2}$, set $d_1^\star \gets d_{ij}^{s_1}$, $d_2^\star \gets w_{i} + w_{j} - d_1^\star + \varepsilon$; otherwise, set $d_2^\star \gets d_{ji}^{s_2}$, $d_1^\star \gets w_{i} + w_{j} - d_2^\star + \varepsilon$; here, $\varepsilon$ is a small positive number. Set $i^\star \gets i$, $j^\star \gets j$.~Add path elimination constraints:
\textit{(i)} $\mathcal{F}^\text{L} \gets \mathcal{F}^\text{L} \cup \left\{(i^\star, j^\star, d_1^\star) \right\}$,
\textit{(ii)} $\mathcal{F}^\text{R} \gets \mathcal{F}^\text{R} \cup \left\{(j^\star, i^\star, d_2^\star) \right\}$.
\end{enumerate}

Set $\mathcal{N}\gets \mathcal{N} \cup \left\{ \left(\tau^\text{L}, \mathcal{F}^\text{L} \right), \left(\tau^\text{R}, \mathcal{F}^\text{R}\right) \right\}$, and go to Step~\ref{main_algo:convergence}.
\end{enumerate}

The modified algorithm gives preference to the branching Step~6$^\prime$c over Steps~\ref{main_algo:branch:continuous} and~\ref{main_algo:branch:discrete}, because unlike the latter, the former would necessarily cut off the VRPTW route set of the parent node in at least one scenario (in both children nodes).
However, as discussed earlier, the path-based disjunctions~\eqref{eq:path_disjunctions} (upon which branching rule~6$^\prime$c is based) are only necessary but not sufficient. This is in contrast to the time window-based disjunctions~\eqref{eq:valid_disjunctions} (upon which the branching rules~\ref{main_algo:branch:continuous} and~\ref{main_algo:branch:discrete} are based), which are both necessary and sufficient. 

\begin{figure}[!htb]
\centering
\begin{tikzpicture}[scale=1,%
					tw/.style = {draw=red!50,text=red},%
					route/.style = {blue!50},%
					projection/.style = {draw=none,ultra thin},%
					customer/.style = {text=blue,circle, inner sep = 0pt, minimum size=0.35cm, fill=white}]
\small
\begin{axis}[name=root,
			 title={$\mathcal{F} = \emptyset$, obj = 25.5},
			 xlabel={Branch as per~\eqref{eq:path_disjunctions}: $(i^\star, j^\star, d_1^\star, d_2^\star) = (4, 3, 7, 5)$}]
	\node[customer] (1s1) at (axis cs:2,10) {1}; \draw[projection] (1s1) -- (axis cs:2,0);
	\node[customer] (2s1) at (axis cs:18,8) {2}; \draw[projection] (2s1) -- (axis cs:18,0);
	\node[customer] (3s1) at (axis cs:14,10){3}; \draw[projection] (3s1) -- (axis cs:14,0);
	\node[customer] (4s1) at (axis cs:7,10) {4}; \draw[projection] (4s1) -- (axis cs:7,0);
	\coordinate (route1s1start) at (axis cs:1,10);
	\coordinate (route2s1start) at (axis cs:1,8);
	\coordinate (route1s1end) at (axis cs:24,10);
	\coordinate (route2s1end) at (axis cs:24,8);
	\draw[route] (axis cs:0,9) -- (route1s1start) -- (1s1) -- (4s1) -- (3s1) -- (route1s1end) -- (axis cs:25,9);
	\draw[route] (axis cs:0,9) -- (route2s1start) -- (2s1) -- (route2s1end) -- (axis cs:25,9);
	
	\node[customer] (1s2) at (axis cs:2,4)  {1}; %\draw[projection] (1s2) -- (axis cs:2,0);
	\node[customer] (2s2) at (axis cs:18,2) {2}; %\draw[projection] (2s2) -- (axis cs:18,0);
	\node[customer] (3s2) at (axis cs:5,2)  {3}; \draw[projection] (3s2) -- (axis cs:5,0);
	\node[customer] (4s2) at (axis cs:12,2) {4}; \draw[projection] (4s2) -- (axis cs:12,0);
	\coordinate (route1s2start) at (axis cs:1,4);
	\coordinate (route2s2start) at (axis cs:1,2);
	\coordinate (route1s2end) at (axis cs:24,4);
	\coordinate (route2s2end) at (axis cs:24,2);
	\draw[route] (axis cs:0,3) -- (route1s2start) -- (1s2) -- (route1s2end) -- (axis cs:25,3);
	\draw[route] (axis cs:0,3) -- (route2s2start) -- (3s2) -- (4s2) -- (2s2) -- (route2s2end) -- (axis cs:25,3);
\end{axis}
\path (root.below south) ++(-4cm,-1.5cm) coordinate (pos1);
\path (root.below south) ++( 4cm,-1.5cm) coordinate (pos2);
\draw[->] (root.below south) -- node[align=left,left]  {\shortstack[l]{No $4-3$ paths with\\travel time $\geq 7$}} +(-2cm,-1.4cm); %to [out=-90,in=0]
\draw[->] (root.below south) -- node[align=right,right] {\shortstack[r]{No $3-4$ paths with\\travel time $\geq 5$}} +(+2cm,-1.4cm); %to [out=-90,in=180]
\begin{axis}[name=l1left,
			 at={(pos1)},
			 anchor=outer north,
 			 title={$\mathcal{F} = \left\{(4, 3, 7)\right\}$, obj = 26.0},
			 xlabel={\bf(optimal)}]
	\node[customer] (1s1) at (axis cs:10,10){1}; \draw[projection] (1s1) -- (axis cs:10,0);
	\node[customer] (2s1) at (axis cs:20,8) {2}; \draw[projection] (2s1) -- (axis cs:20,0);
	\node[customer] (3s1) at (axis cs:4,10) {3}; \draw[projection] (3s1) -- (axis cs:4,0);
	\node[customer] (4s1) at (axis cs:15,10){4}; \draw[projection] (4s1) -- (axis cs:15,0);
	\coordinate (route1s1start) at (axis cs:1,10);
	\coordinate (route2s1start) at (axis cs:1,8);
	\coordinate (route1s1end) at (axis cs:24,10);
	\coordinate (route2s1end) at (axis cs:24,8);
	\draw[route] (axis cs:0,9) -- (route1s1start) -- (3s1) -- (1s1) -- (4s1) -- (route1s1end) -- (axis cs:25,9);
	\draw[route] (axis cs:0,9) -- (route2s1start) -- (2s1) -- (route2s1end) -- (axis cs:25,9);
	
	\node[customer] (1s2) at (axis cs:10,4) {1}; %\draw[projection] (1s2) -- (axis cs:10,0);
	\node[customer] (2s2) at (axis cs:20,2) {2}; %\draw[projection] (2s2) -- (axis cs:20,0);
	\node[customer] (3s2) at (axis cs:4,2)  {3}; %\draw[projection] (3s2) -- (axis cs:4,0);
	\node[customer] (4s2) at (axis cs:14,2) {4}; \draw[projection] (4s2) -- (axis cs:14,0);
	\coordinate (route1s2start) at (axis cs:1,4);
	\coordinate (route2s2start) at (axis cs:1,2);
	\coordinate (route1s2end) at (axis cs:24,4);
	\coordinate (route2s2end) at (axis cs:24,2);
	\draw[route] (axis cs:0,3) -- (route1s2start) -- (1s2) -- (route1s2end) -- (axis cs:25,3);
	\draw[route] (axis cs:0,3) -- (route2s2start) -- (3s2) -- (4s2) -- (2s2) -- (route2s2end) -- (axis cs:25,3);

	\draw[tw] (axis cs:4,5)  rectangle node {$\tau_3^\star$} +(axis cs:3,2);
	\draw[tw] (axis cs:10,5) rectangle node {$\tau_1^\star$} +(axis cs:3,2);
	\draw[tw,draw=brown] (axis cs:14,5) rectangle +(axis cs:8,2);
	\node[text=brown,right] at (axis cs:21.5,6) {$\leftarrow \tau_4^\star$};
	\draw[tw] (axis cs:17,5) rectangle node {$\tau_2^\star$} +(axis cs:3,2);
\end{axis}
\begin{axis}[name=l1right,
			 at={(pos2)},
			 anchor=outer north,
			 title={$\mathcal{F} = \left\{(3, 4, 5)\right\}$, obj = 26.0},
			 xlabel={\bf(also optimal)},
			 ytick={2.5,9}]
	\node[customer] (1s1) at (axis cs:2,10) {1}; \draw[projection] (1s1) -- (axis cs:2,0);
	\node[customer] (2s1) at (axis cs:18,8) {2}; \draw[projection] (2s1) -- (axis cs:18,0);
	\node[customer] (3s1) at (axis cs:14,10){3}; \draw[projection] (3s1) -- (axis cs:14,0);
	\node[customer] (4s1) at (axis cs:7,10) {4}; \draw[projection] (4s1) -- (axis cs:7,0);
	\coordinate (route1s1start) at (axis cs:1,10);
	\coordinate (route2s1start) at (axis cs:1,8);
	\coordinate (route1s1end) at (axis cs:24,10);
	\coordinate (route2s1end) at (axis cs:24,8);
	\draw[route] (axis cs:0,9) -- (route1s1start) -- (1s1) -- (4s1) -- (3s1) -- (route1s1end) -- (axis cs:25,9);
	\draw[route] (axis cs:0,9) -- (route2s1start) -- (2s1) -- (route2s1end) -- (axis cs:25,9);
	
	\node[customer] (1s2) at (axis cs:2,4)  {1};
	\node[customer] (2s2) at (axis cs:18,2.5){2};
	\node[customer] (3s2) at (axis cs:14,1) {3}; 
	\node[customer] (4s2) at (axis cs:7,2.5){4}; 
	\coordinate (route1s2start) at (axis cs:1,4);
	\coordinate (route2s2start) at (axis cs:1,2.5);
	\coordinate (route3s2start) at (axis cs:1,1);
	\coordinate (route1s2end) at (axis cs:24,4);
	\coordinate (route2s2end) at (axis cs:24,2.5);
	\coordinate (route3s2end) at (axis cs:24,1);
	\draw[route] (axis cs:0,2.5) -- (route1s2start) -- (1s2) -- (route1s2end) -- (axis cs:25,2.5);
	\draw[route] (axis cs:0,2.5) -- (route2s2start) -- (4s2) -- (2s2) -- (route2s2end) -- (axis cs:25,2.5);
	\draw[route] (axis cs:0,2.5) -- (route3s2start) -- (3s2) -- (route3s2end) -- (axis cs:25,2.5);
	
	\draw[tw,draw=brown!70] (axis cs:0,5)  rectangle +(axis cs:7,2);
	\node[text=brown,right] at (axis cs:6.5,6) {$\leftarrow \tau_4^\star$};
	\draw[tw] (axis cs:12,5) rectangle node {$\tau_3^\star$} +(axis cs:3,2);
	\draw[tw] (axis cs:2,5)  rectangle node {$\tau_1^\star$} +(axis cs:3,2);
	\draw[tw] (axis cs:17,5) rectangle node {$\tau_2^\star$} +(axis cs:3,2);
\end{axis}
\end{tikzpicture}
\caption{The search tree of our algorithm (utilizing path-based disjunctions) for the illustrative example of Figure~\ref{figure:example_instance}.}
\label{figure:example_algorithm_with_path_disjunctions}
\end{figure}
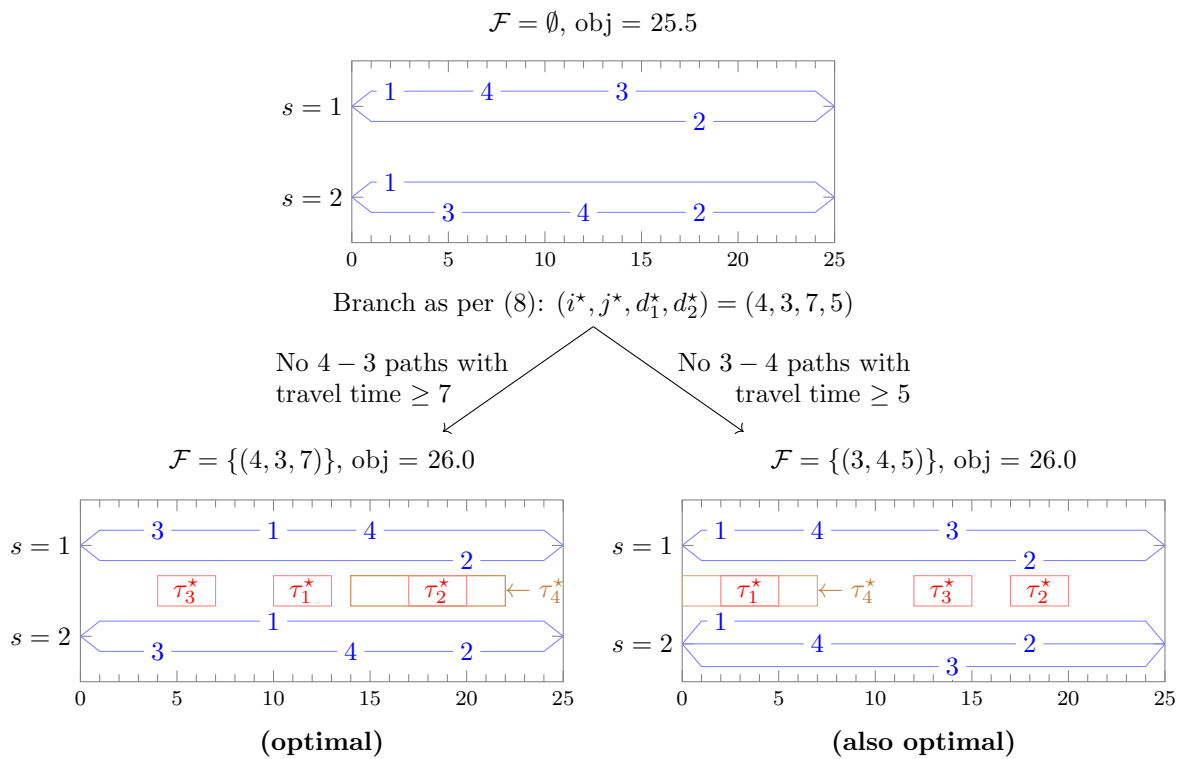

Figure~\ref{figure:example_algorithm_with_path_disjunctions} shows the effect that the new branching rules have on the branch-and-bound search tree for the case of our illustrative example from Figure~\ref{figure:example_instance}.
Note how, in the root node, instead of branching via the time window-based disjunctions~\eqref{eq:valid_disjunctions}, the modified branching Step~\ref{path_based_algo:branch} certifies that it is impossible to have both $3-4$ and $4-3$ paths in different scenario route sets, and branches using the disjunction~\eqref{eq:path_disjunctions} instead. Observe that the search tree is much smaller than in the case of Figure~\ref{figure:example_algorithm}. Our numerical experiments confirm that this is generally true, i.e., that incorporating the path-based disjunctions~\eqref{eq:path_disjunctions} results in fewer nodes being explored (see Section~\ref{sec:computational_results_detailed}).

\subsection{Generating Upper Bounds}\label{sec:solution_approach_upper_bounds}
A small value of the global upper bound $UB$ can significantly speed up the solution process by fathoming more nodes of the search tree using the fathoming Step~\ref{main_algo:fathom_by_bound}.
The algorithm presented in Section~\ref{sec:solution_approach_algorithm} updates $UB$ in Step~\ref{main_algo:check_feasibility} only.
However, it is possible to update $UB$ more frequently by generating candidate (feasible) TWAVRP solutions using the route sets $\{\mathbf{R}_s\}_{s \in \mathcal{S}}$ obtained in the node processing Step~\ref{main_algo:node_processing}.
The basic idea is to select some scenario and use its solution as a ``template'' \cite{Groer2009}, assigning the time windows based on the arrival times in this template solution.
Specifically, for every scenario $s^\star \in \mathcal{S}$ in which a new route set $\mathbf{R}_{s^\star}$ is computed in Step~\ref{main_algo:node_processing} of the algorithm, we can attempt to generate a feasible TWAVRP solution $(\tau^h, \{\mathbf{R}_s^h\}_{s \in \mathcal{S}})$ using the following procedure.
\begin{enumerate}
\item\label{generating_upper_bounds:init}
\textit{Initialize arrival times.} Let $\mathbf{R} \gets \mathbf{R}_{s^\star}$ and $\theta \gets \theta_{s^\star}$. Let $\bm{a} \in \mathcal{X}(\mathbf{R}, [e, \ell]; \theta)$ be the vector of arrival times with minimum cumulative waiting time. That is, we recursively define $a_{R_{k,0}} = e_0$ and $a_{R_{k,l+1}} = \max\{e_{R_{k,l+1}}, a_{R_{k,l}} + t_{R_{k,l},R_{k,l+1}}\}$ for all $l \in \{1, \ldots, n_k - 1\}$ and $k \in \{1, \ldots, m\}$. For $i \in V_C$ such that $q_i = 0$ (that is, $i \neq R_{k,l}$ for any $k,l$), we define $a_i = e_i$.

\item\label{generating_upper_bounds:assign_tw}
\textit{Assign time windows.} Let $\tau^h$ be defined as follows.
\begin{equation*}
\tau_i^h = \begin{dcases}
\left[x_i, x_i + w_i\right],
\text{where } x_i = \begin{dcases}
e_i          & \text{if } a_i - w_i/2 \leq e_i \\
\ell_i - w_i & \text{if } a_i + w_i/2 \geq \ell_i \\
a_i - w_i/2  & \text{otherwise}
\end{dcases} & \text{if } i \in V_\text{cont} \\
\left[x_{i,b_i}, y_{i,b_i}\right], \text{where } b_i \in \mathop{\arg\min}_{b \in \{1,\ldots,N_i\}} \left\{ \min_{\omega \in [x_{i,b}, y_{i,b}]} \left\lvert a_i - \omega \right\rvert \right\} & \text{if } i \in V_\text{disc}
\end{dcases} \quad \forall i \in V_C.
\end{equation*}

\item\label{generating_upper_bounds:compute_ub}
\textit{Compute upper bound.}
For each $s \in \mathcal{S}$, let $\mathbf{R}^h_s$ be a (possibly suboptimal) solution of $\mathtt{VRPTW}(\tau^h; \theta_{s})$. Let $ub^h \gets \sum_{s \in \mathcal{S}} p_{s} c(\mathbf{R}^h_{s})$. If $ub^h \leq UB$, set $UB \gets ub^h$ and $\tau^\star \gets \tau^h$.
\end{enumerate}
We remark that it is not necessary to solve the VRPTW instances to optimality in Step~\ref{generating_upper_bounds:compute_ub} above, since the generated upper bounds $ub^h$ are guaranteed to still be valid. Consequently, we can utilize any (possibly heuristic) VRPTW solver to quickly compute candidate upper bounds.
For example, the computational results reported in Section~\ref{sec:computational_results} were obtained by implementing Solomon's sequential insertion construction heuristic ``I1''~\cite{Solomon1987} in combination with a local search procedure~\cite{Braysay2005}, which used the \textit{Relocate}, \textit{2-opt}, \textit{2-opt$^\star$} and \textit{Or-opt} moves within a deterministic Variable Neighborhood Descent algorithm (e.g., see~\cite{Braysay2003}).

\subsection{Modification as a Heuristic Algorithm}\label{sec:solution_approach_heuristic_modification}
The exact algorithm of Section~\ref{sec:solution_approach_algorithm} can be readily modified as a heuristic algorithm. Indeed, if one uses a heuristic VRPTW solver (e.g., one based on a metaheuristic) in place of an exact one in the node processing Step~\ref{main_algo:node_processing}, then the time window assignment $\tau^\star$ determined by the algorithm is still guaranteed to be feasible, although not necessarily optimal.
Overall TWAVRP optimality can no longer be guaranteed, since the fathoming Step~\ref{main_algo:fathom_by_bound} may incorrectly prune a node whose descendant contains the optimal time window assignment.
Nevertheless, this modification as a heuristic is particularly suited from a practical viewpoint, as it allows one to utilize any available VRPTW heuristic solver ``out~of~the~box.''

\section{Exact Solution of VRPTW Subproblems}\label{sec:implementation_details}
Various exact solution schemes have been proposed for the VRPTW over the last several decades.
These include algorithms based on branch-and-cut, Lagrangean relaxation and column generation, among others.
We refer the reader to~\cite{chapterInTothVigo} for a recent survey.
The most successful of these are \textit{branch-price-and-cut} algorithms, which correspond to branch-and-bound algorithms in which the bounds are obtained by solving linear relaxations of a \emph{set partitioning} model by column generation, and are further strengthened by generating cutting planes.

In order to solve VRPTW instances in Steps~\ref{main_algo:node_processing} and~\ref{path_based_algo:node_processing} of the main algorithm, we implemented the branch-price-and-cut algorithm described in~\cite{pecin2017improved}, as well as procedures to warm start the latter exact method in the context of our algorithm. Our implementation incorporates several elements of the algorithm described in~\cite{pecin2017improved}, including $ng$-routes, bidirectional labeling, variable fixing, route enumeration, and limited-memory subset row cuts. In what follows, we highlight only the most important of these ingredients; for details, we refer the reader to~\cite{pecin2017improved}.

\subsection{Branch-Price-and-Cut Implementation}\label{sec:implementation_details_vrptw_spp}
For ease of notation, we shall drop the subscript $s$ referencing a scenario and assume that all operational parameters $\theta$ are fixed to certain given values. We shall also assume that the time windows have been fixed at $\tau$ and that the set of forbidden paths is given to be $\mathcal{F}$.
Note that we only describe the solution approach for~\ref{eq:model_vrptw_with_path_elimination}; the approach for~\ref{eq:model_vrptw} is obtained by simply setting $\mathcal{F} = \emptyset$.
We remark that, before we call the exact solution method, we modify the VRPTW instance to obtain an equivalent one by tightening the time windows $\tau$ and reducing the arc set $A$ using the preprocessing routines described in~\cite{Ascheuer2001,Kallehauge2007}.
In addition to this, it is possible to further reduce the arc set $A$ using members of $\mathcal{F}$. Indeed, if for some $(i, j, d) \in \mathcal{F}$, the shortest travel time from $i$ to $j$ in graph $G$ exceeds $d$, then we can remove the arc $(i, j)$ from $A$.

In the following, $\mathcal{P}_{ij}$ denotes the set of all $i-j$ paths in graph $G$ (after preprocessing), where $(i, j) \in V_C \times V_C$, such that $i\neq j$. We use $\Omega$ to denote the set of all (elementary) vehicle routes that are feasible with respect to capacity and time window constraints.
For a given route $r \in \Omega$, $\lambda_{ir}$ denotes the number of times customer $i \in V_C$ is visited in route $r$,  $\eta_{ijr}$ denotes the number of times arc $(i, j) \in A$ is traversed by route $r$, while $c_r$ denotes its cost, i.e., $c_r \equiv c(r)$.
The set partitioning model is described in the following.
In this model, $x_r$ is a binary path-flow variable that encodes whether route $r \in \Omega$ is part of the optimal route set.
\begin{equation}\label{eq:set_partitioning_model}
\begin{array}{r@{\;\;}l}
\displaystyle \mathop{\text{minimize}}_{x} & \displaystyle \sum_{r \in \Omega} c_r x_r \\
\text{subject to} & x_r \in \{0, 1\}, \quad \forall r \in \Omega, \\
& \displaystyle \sum_{r \in \Omega} \lambda_{ir} x_r = 1, \quad \forall i \in V_C, \\
& \displaystyle \sum_{r \in \Omega} \sum_{(i^\prime,j^\prime) \in \pi} \eta_{i^\prime j^\prime r} x_r \leq \lvert \pi \rvert - 1, \quad \forall \pi \in \mathcal{P}_{ij}: t(\pi) \geq d, \;\; \forall (i, j, d) \in \mathcal{F}.
\end{array}
\end{equation}
In the above, the last set of inequalities are \emph{infeasible path elimination constraints} (e.g., see~\cite{Ascheuer2001,Kallehauge2007}) that forbid the occurrence of $i - j$ paths with travel time greater than or equal to $d$.
In our implementation, we replace the subscript of the innermost summation with $(i^\prime, j^\prime) \in \text{tr.cl.}(\pi)$, where $\text{tr.cl.}(\pi)$ denotes the transitive closure of $\pi$.\footnote{If $\pi = (v_1, \ldots, v_p)$ denotes an elementary path, then its transitive closure is the set of arcs $(v_k, v_l)$ such that $v_l$ can be reached from $v_k$ using only arcs in $\pi$, i.e, $\text{tr.cl.}(\pi) = \{(v_k, v_l) \in A : (k, l) \in \{1, \ldots, p\} \times \{1, \ldots, p\}, k < l \}$.} This so-called \emph{tournament form} of the inequality is stronger than the version presented above, see~\cite{Ascheuer2001} for a proof. Furthermore, it is also well known that one can relax the feasible space of the above set partitioning model by including non-elementary vehicle routes in $\Omega$ without sacrificing optimality. In our implementation, we replace $\Omega$ with the set of so-called $ng$-routes $\Omega^\text{ng} \supseteq \Omega$, which are not necessarily elementary~\cite{baldacci2011new}.

We now describe the branch-price-and-cut algorithm to solve the set partitioning model over $\Omega^\text{ng}$.
The root node, which solves the linear relaxation of this set partitioning model, is initialized with a subset of $\Omega^\text{ng}$ (single-customer vehicle routes) but no infeasible path elimination constraints.
A pricing subproblem is used to generate other members of $\Omega^\text{ng}$ (also referred to as columns), as necessary.
After column generation has converged, if the gap between the current node lower bound and global upper bound is sufficiently small ($\leq 1\%$ in our implementation), we employ \textit{route enumeration} to generate all feasible vehicle routes with reduced costs smaller than this gap~\cite{baldacci2008exact}. \cite{baldacci2008exact} have shown that this subset must contain the routes of all optimal solutions. Therefore, we can solve the resulting ``reduced'' set partitioning model by using a standard integer programming solver.
We remark that this is done only if the number of generated routes is less than a threshold, which we set to $3\times 10^6$, as suggested in~\cite{pecin2017improved}.
On the other hand, if the current node gap is large or if route enumeration generated too many routes, then we attempt to separate the infeasible path elimination constraints (as well as other valid inequalities) in order to tighten the linear relaxation.
The above process is iterated until we cannot generate any more columns or inequalities.
At this stage, if the current node solution is fractional, we create additional children nodes by branching on the number of used vehicles or by branching on edges/arcs.

\paragraph{Pricing subproblem.} The pricing subproblem is a \textit{shortest path problem with resource constraints}, with customer demands and arc travel times considered as resources constrained by the vehicle capacity and time windows, respectively.
We utilize the dynamic programming algorithm described in~\cite{pecin2017improved} to solve this problem.
In order to speed up the solution of the pricing subproblem, we apply various techniques including bidirectional labeling and variable fixing (based on reduced costs). In addition, we implemented the bucket pruning heuristic (e.g., see~\cite{Fukasawa2006}) to find candidate columns and use the dynamic programming algorithm only if the former fails to generate columns.

\paragraph{Route enumeration.} We use the dynamic programming algorithm of~\cite{pecin2017improved} with a modification to account for the presence of the infeasible path elimination constraints. In particular, we consider two different ``partial routes'' that visit the same set of customers (but possibly in different sequences) to be undominated irrespectively of their resource consumptions, and we do not perform any associated dominance checks.
This prevents incorrectly ``pruning'' a vehicle route that satisfies the infeasible path elimination constraints in favor of one that does not.
Among all enumerated routes, we only consider those which are elementary and satisfy all infeasible path elimination constraints (which can be done in $\mathcal{O}(\left\lvert \mathcal{F} \right\rvert)$ time per route) to include in the final ``reduced'' set partitioning model.

\paragraph{Cutting planes.} We use the tournament form of the infeasible path elimination constraints as ``necessary cuts'' (a violating member is guaranteed to be separated), and the \textit{extended capacity cuts}~\cite{pessoa2008robust} and \textit{limited-node-memory subset row cuts}~\cite{pecin2017improved} as ``strengthening cuts'' (a violating member may not necessarily be separated). The separation algorithms for these cuts are described next.
\begin{itemize}
\item The infeasible path elimination constraints are separated by utilizing the polynomial-time path-growing scheme of~\cite{Ascheuer2001}.
Specifically, suppose $G_{\bar{x}} = (V, A_{\bar{x}})$ is the so-called support graph of the current linear programming solution $\bar{x}$, where $A_{\bar{x}} = \left\{ (i^\prime, j^\prime) \in A: \sum_{r \in \Omega^\text{ng}} \eta_{i^\prime j^\prime r} \bar{x}_r > 0 \right\}$.
Then, for every $(i, j, d) \in \mathcal{F}$, the scheme of~\cite{Ascheuer2001} is used to obtain the set of all $i-j$ paths $\pi$ in $G_{\bar{x}}$, which satisfy $\sum_{r \in \Omega^\text{ng}} \sum_{(i^\prime,j^\prime) \in \text{tr.cl.}(\pi)} \eta_{i^\prime j^\prime r} \bar{x}_r > \lvert \pi \rvert - 1$.
Amongst all such paths, we choose the ones for which the travel time $t(\pi)$ is greater than $d$ and add the corresponding tournament form of the infeasible path elimination constraints to the current linear relaxation.

\item The extended capacity cuts are separated by first using the CVRPSEP package~\cite{Lysgaard2004} to separate the so-called \textit{rounded capacity inequalities} and then lifting these, as described in~\cite{pessoa2008robust}.

\item The limited-node-memory subset row cuts are separated by first identifying all subset row cuts violated by node sets with cardinality up to $5$ and then identifying their so-called ``node memory sets,'' as described in~\cite{pecin2017improved}.
\end{itemize}
We remark that, since the infeasible path elimination inequalities are defined over arcs, they are ``robust'' and affect the pricing subproblem only through a corresponding term in the arc cost, i.e., their dual value. In other words, the addition of these inequalities does not affect the complexity of the pricing subproblem.
In contrast, the extended capacity and limited-node-memory subset row cuts are ``non-robust'' as their addition increases the complexity of the pricing subproblem~\cite{pecin2017improved}.

\subsection{Warm Starting}\label{sec:implementation_details_vrptw_warmstarting}
Initializing the branch-price-and-cut algorithm described in the previous section with a feasible set of columns can speed up the convergence of its column generation process, leading to small computation times.
In addition to this, providing a valid initial upper bound on the optimal objective value can also significantly speed up the search, both in the context of route enumeration, where it can result in fewer routes being enumerated, as well as branch-and-bound, where more parts of the search tree can be fathomed early in the process.
In this section, we describe warm starting procedures that can be employed in the context of our algorithm of Section~\ref{sec:solution_approach}.

Consider a parent node $(\tau, \mathcal{F})$ (already processed) with optimal solution $\{\mathbf{R}_s\}_{s\in \mathcal{S}}$ as obtained in Step~\ref{main_algo:node_processing}.
Also, consider one of its child nodes $(\tau^h, \mathcal{F}^h)$, where $h \in \{\text{L}, \text{R}\}$, created in the branching Step~\ref{path_based_algo:branch}.
Finally, consider a scenario $s \in \mathcal{S}$ such that the corresponding optimal route set in the parent node, $\mathbf{R}_s$, is not feasible for the child node $(\tau^h, \mathcal{F}^h)$.
This means that $\mathbf{R}_s$ cannot be directly transferred to the latter, and warm starting is sought to aid the search towards the optimal one.

\paragraph{Generating an initial set of columns.} Let $\Omega^\text{ng}_s$ be the set of columns that were generated during the branch-price-and-cut process in scenario $s$ of the parent node $(\tau, \mathcal{F})$.
Since the child node $(\tau^h, \mathcal{F}^h)$ differs from its parent in exactly one constraint, it is likely that several members of $\Omega^\text{ng}_s$ are also feasible in the set partitioning model of the child's scenario$-s$ VRPTW subproblem.
Therefore, we can simply loop through the members of $\Omega^\text{ng}_s$ and filter out all infeasible columns to generate the initial linear relaxation in the branch-price-and-cut algorithm.

\paragraph{Generating an initial upper bound.} The following procedure computes a valid upper bound, $ub_s^h$, on the cost of the child's scenario$-s$ VRPTW subproblem. The procedure attempts to ``repair'' the scenario$-s$ optimal route set of the parent node and generate one that is feasible for the child.
\begin{enumerate}
\item\label{warm_start:init}
Set $ub_s^h \gets UB - \sum_{s^\prime \in \mathcal{S}\setminus \{s\}} c(\mathbf{R}_{s^\prime})$, where $UB$ is the currently applicable upper bound from the algorithm of Section~\ref{sec:solution_approach}; and set $\mathbf{R}^\prime \gets \mathbf{R}_s$.

\item\label{warm_start:if}
Let $H \subseteq V_C$ be defined as follows: if $(\tau^h, \mathcal{F}^h)$ was created using branching Steps~\ref{main_algo:branch:continuous} or~\ref{main_algo:branch:discrete}, then $H = \{i^\star\}$; otherwise $H = \{i^\star, j^\star\}$.

\item\label{warm_start:local_search}
For each $i \in H$:
\begin{enumerate}
\item Remove customer $i$ from its current position in $\mathbf{R}^\prime$ and insert it into a new vehicle route.

\item Apply local search on $\mathbf{R}^\prime$, ensuring that each accepted move satisfies all time window and path elimination constraints in $(\tau^h, \mathcal{F}^h)$.

\item If $c(\mathbf{R}^\prime) < ub_s^h$, set $ub_s^h \gets c(\mathbf{R}^\prime)$.
\end{enumerate}
\end{enumerate}
In our implementation of local search, we considered the \textit{Relocate}, \textit{2-opt}, \textit{2-opt$^\star$} and \textit{Or-opt} moves within a deterministic Variable Neighborhood Descent algorithm (e.g., see~\cite{Braysay2003}).
We remark that the solution of (up to $S$) VRPTW instances in Step~\ref{main_algo:node_processing} of the algorithm can be easily parallelized.
Alternatively, one can do this serially on a single CPU thread (e.g., by starting with the lowest indexed scenario). 
In the former setting, no information can be exchanged among the VRPTW instances.
In the latter setting, however, one can capitalize on instances that have already been solved so as to obtain improved upper bounds in Step~\ref{warm_start:init} of the above procedure as follows:
\[
ub_s^h \gets UB - \sum_{s^\prime \in \mathcal{S}: s^\prime < s} c(\mathbf{R}_{s^\prime}^h) - \sum_{s^\prime \in \mathcal{S}: s^\prime > s} c(\mathbf{R}_{s^\prime}),
\]
where $\mathbf{R}_{s^\prime}^h$ is the just computed optimal solution in scenario $s^\prime$ of the node under consideration.

\section{Computational Results}\label{sec:computational_results}
This section presents computational results obtained by our algorithm on benchmark instances from the literature. Specifically, in Section~\ref{sec:computational_results_instances}, we present the characteristics of the test instances; in Section~\ref{sec:computational_results_comparison}, we present a summary of the numerical performance of our algorithm and compare it to existing solution methods; in Section~\ref{sec:computational_results_detailed}, we present detailed tables of results outlining the performance of each component of our algorithm; and, finally in Section~\ref{sec:computational_results_large_instances}, we present results of a parallel implementation on instances containing a large number of scenarios.

The algorithm was coded in C++ and the runs were conducted on an Intel Xeon E5-2687W 3.1~GHz processor with 4~GB of available RAM. 
The nodes in Step~\ref{main_algo:convergence} of the algorithm were selected using a simple depth-first rule that backtracked whenever the gap between the objective value of the current node and the current upper bound exceeded 50\% of the gap between the global lower and current upper bound.
All subordinate linear and mixed-integer linear programs were solved using default settings of the IBM ILOG CPLEX Optimizer~12.7.
Finally, except for the results presented in Section~\ref{sec:computational_results_large_instances}, all runs were restricted to a single CPU thread. This facilitates a fair comparison with existing algorithms from the literature in Section~\ref{sec:computational_results_comparison} and between different configurations of our algorithm in Section~\ref{sec:computational_results_detailed}. The results presented in Section~\ref{sec:computational_results_large_instances} were obtained with OpenMP by using up to $\min\{S, 10\}$ threads in parallel, where $S$ is the number of scenarios. In all cases, an overall ``wall clock'' time limit of one hour per instance was imposed.

\subsection{Benchmark Instances}\label{sec:computational_results_instances}
Existing benchmark instances for the TWAVRP focus solely on demand uncertainty. 
For the continuous setting, the authors of~\cite{Spliet2015_continuous} introduced 40 randomly generated instances. 
The number of customers ($n$) in these instances varies from 10 to 25. Subsequently, \cite{Dalmeijer2018} proposed 50 additional instances with $n$ varying from 30 to 50. Each of the 90 available instances consists of 3 demand scenarios (low, medium, high), each with equal probability of occurrence. The average demand (for each customer) across the three scenarios is about $1/6$ of the vehicle capacity $Q$.
The exogenous time windows are designed to be much wider than the endogenous time windows; in particular, the average (across customers) exogenous time window width ($\ell_i - e_i$) is 10.8, compared to an endogenous time window width ($w_i$) of just 2.0.
For the discrete setting, the authors of~\cite{Spliet2015_discrete} introduced 80 randomly generated instances with $n$ varying from 10 to 60. Except for the structure of the feasible time window set, the instances share similar characteristics as in the continuous setting, each consisting of 3 demand scenarios. In each instance, the number of candidate time windows ($N_i$) is equal to 3 for about 10\% of the customers, 5 for about 60\% of the customers, and 7 for the remaining 30\%.
Similarly to the continuous setting, the customer locations are generated using a uniform distribution over a square with the depot located in the center.
Moreover, the time windows and vehicle capacities are chosen such that no more than eight customers can be visited on a single vehicle route in any scenario.
All of the aforementioned test instances are inspired from the Dutch retail sector, and can be found online at {http://people.few.eur.nl/spliet}.

To test our algorithm on instances containing a large number of scenarios, we generated 80 additional benchmark instances.
Specifically, for each existing (continuous and discrete) TWAVRP instance with $n \leq 25$, we used a similar procedure as described in~\cite{Spliet2015_discrete} to generate 15 additional demand scenarios.
For a particular instance, we first generate a nominal demand $\bar{q}_i$ for each customer $i \in V_C$ using a normal distribution with mean $5.0$ and variance $1.5$.
To generate additional scenarios, we draw additive disturbances $\epsilon_{si}$ from a uniform distribution on $[-1.5, 1.5]$ for each $i \in V_C$ and $s \in \mathcal{S}$, and multiplicative factors $f_s$ from a uniform distribution on $[0.625, 1.375]$ for each $s \in \mathcal{S}$.
The demand of customer $i$ in scenario $s$ is then computed as $q_{si} = \left\lceil \max\left\{ f_s \left(\bar{q}_i + \epsilon_{si} \right), 10^{-6} \right\} \right\rceil$.
The multiplicative factors $f_s$ determine the level of correlation among the customer demands.
For instance, high (low) values of $f_s$ may represent the behavior that demands increase (decrease) uniformly for all retailers in a supply chain, whereas values close to one represent the nominal situation in which the demands are uncorrelated.
The additional benchmark instances can be downloaded from {http://gounaris.cheme.cmu.edu/datasets/twavrp}.

\subsection{Comparison with Existing Methods}\label{sec:computational_results_comparison}
\begin{figure}[!htb]
\captionsetup[subfigure]{belowskip=0pt}
\centering
\begin{subfigure}[b]{0.5\textwidth}\centering
\includegraphics[height=6.2cm]{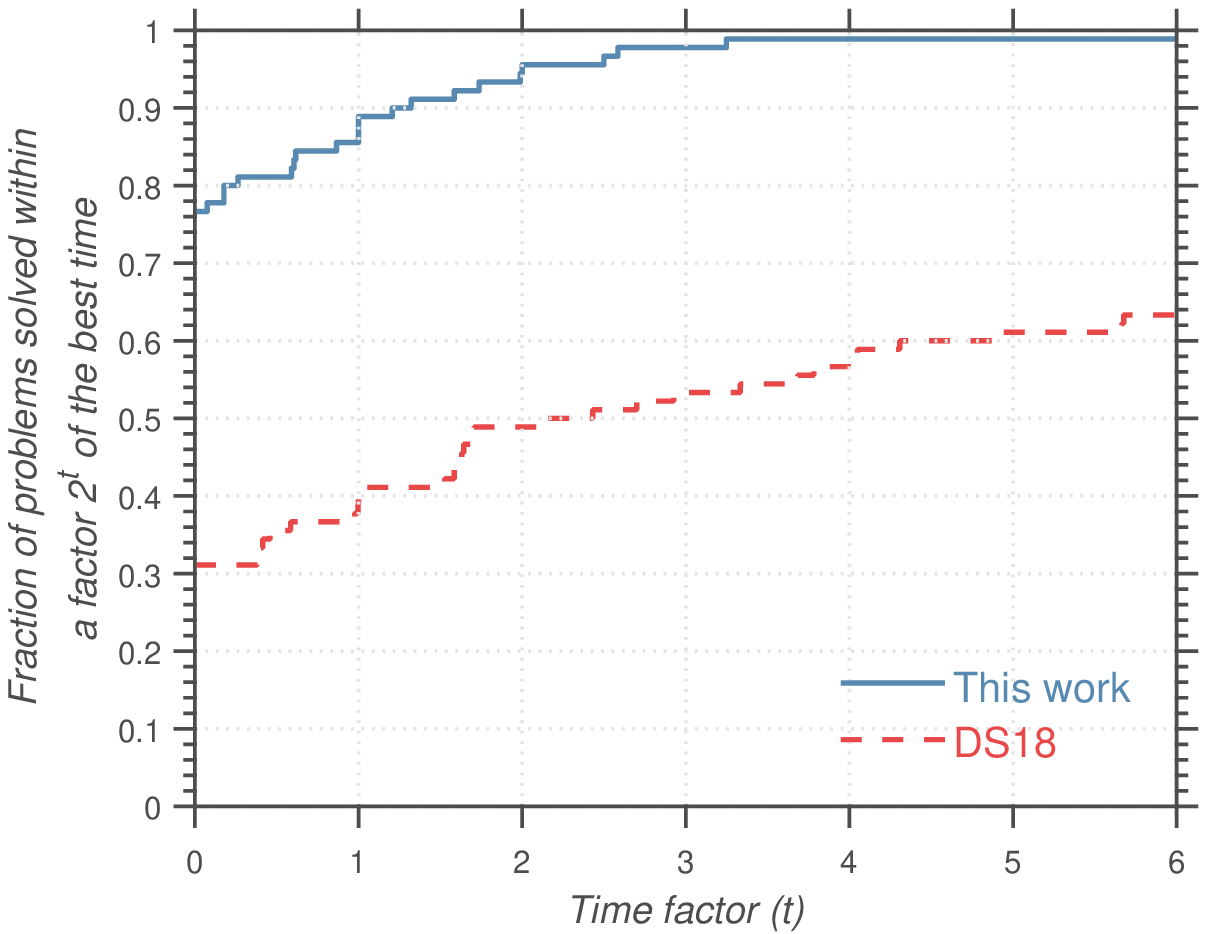}
\caption{Continuous}
\label{figure:perf_compare_with_existing_continuous}
\end{subfigure}~%
\begin{subfigure}[b]{0.5\textwidth}\centering
\includegraphics[height=6.2cm]{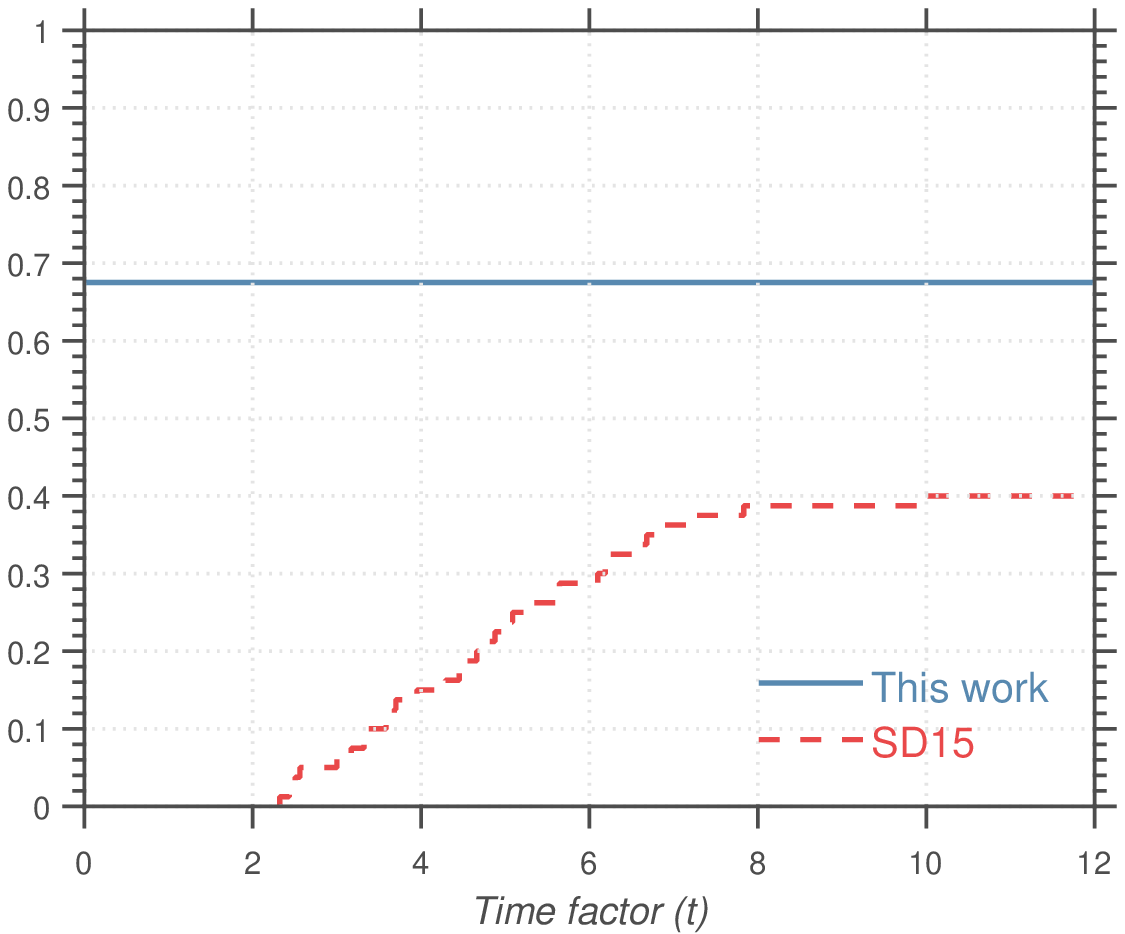}
\caption{Discrete}
\label{figure:perf_compare_with_existing_discrete}
\end{subfigure}
\caption{Log-scaled performance profiles across all benchmark instances. The left graph compares the performance in the continuous setting, in which ``DS18'' shows the performance of the algorithm of~\cite{Dalmeijer2018}, while the right graph compares the performance in the discrete setting, in which ``SD15'' shows the performance of the algorithm of~\cite{Spliet2015_discrete}. In both graphs, ``This work'' shows the performance of our proposed algorithm. For each curve (i.e., algorithm), the value at $t = 0$ gives the fraction of benchmark instances for which it is fastest, while the limiting value at $t \to \infty$ gives the fraction of instances which it could solve within the time limit of one hour.}
\label{figure:perf_compare_with_existing}
\end{figure}

We first compare the performance of our algorithm with the results published in~\cite{Dalmeijer2018} for the case of the continuous TWAVRP. We do not compare with the algorithm of~\cite{Spliet2015_continuous}, since the authors of~\cite{Dalmeijer2018} have demonstrated that their algorithm is superior to the former. Table~\ref{table:computational_comparison_continuous} summarizes the comparison of the numerical performance across all 90 instances that are available for the continuous TWAVRP.
The column \textbf{\#} denotes the number of test instances that contain $n$ customers. For each algorithm, \textbf{Optimal} denotes the number of test instances that it could solve to optimality in one hour while \textbf{Time (sec)} denotes the average time in seconds to solve these instances to optimality.
For those instances which could not be solved to optimality in one hour, the column \textbf{Gap (\%)} reports the average optimality gap, defined as $(UB - LB)/UB \times 100\%$, where $LB$ and $UB$ are respectively the global lower and upper bounds determined by the algorithm after one hour.
The two methods are also compared in the performance profiles~\cite{Dolan2002} of Figure~\ref{figure:perf_compare_with_existing_continuous}.
Our proposed algorithm is able to solve all but one (89 out of 90) benchmark instances to optimality, utilizing an average computation time of 169 seconds; of these, 32 instances were unsolved by the best previous method, while the one unsolved instance was determined to be within 0.8\% of optimality.
These results demonstrate that our algorithm strongly outperforms the existing method, solving more problems and achieving (or matching) the fastest computation time in all instances.

% Table generated by Excel2LaTeX from sheet 'Sheet2'
\begin{table}[!htb]
  \centering
  \caption{Computational comparison of the proposed algorithm against the existing state-of-the-art algorithm (DS18)~\cite{Dalmeijer2018} on all 90 benchmark instances of the continuous TWAVRP.}
    \begin{tabularx}{\textwidth}{lc*{6}{R}}
    \toprule
           &        & \multicolumn{3}{c}{DS18} & \multicolumn{3}{c}{Proposed algorithm} \\
    \cmidrule(lr){3-5}\cmidrule(l){6-8}
    $n$    & \#     & Optimal & Time (sec) & Gap (\%) & Optimal & Time (sec) & Gap (\%) \\
    \midrule
    10     & 10     & 10     & 0.1    & --     & 10     & 0.1    & -- \\
    15     & 10     & 10     & 4.6    & --     & 10     & 0.6    & -- \\
    20     & 10     & 10     & 2.2    & --     & 10     & 1.5    & -- \\
    25     & 10     & 10     & 12.4   & --     & 10     & 8.6    & -- \\
    30     & 10     & 9      & 204.5  & 1.66   & 10     & 48.2   & -- \\
    35     & 10     & 6      & 152.8  & 0.89   & 9      & 51.5   & 0.79 \\
    40     & 10     & 2      & 1860.0 & 1.16   & 10     & 342.3  & -- \\
    45     & 10     & 0      & --     & 2.74   & 10     & 361.3  & -- \\
    50     & 10     & 0      & --     & 4.26   & 10     & 696.0  & -- \\[0.1cm]
    All    & 90     & 57     & 117.0  & 2.56   & 89     & 169.1  & 0.79 \\
    Processor &        & \multicolumn{3}{r}{Intel i7 3.5 GHz} & \multicolumn{3}{r}{Xeon E5-2687W 3.1GHz} \\
    \bottomrule
    \end{tabularx}%
  \label{table:computational_comparison_continuous}%
\end{table}%

We now turn our attention to the discrete TWAVRP. Table~\ref{table:computational_comparison_discrete} compares the numerical performance of our algorithm with the results published in~\cite{Spliet2015_discrete} across all 80 instances of the discrete TWAVRP. The columns in this table have the same meaning as in Table~\ref{table:computational_comparison_continuous}. The two algorithms are also compared in the performance profiles of Figure~\ref{figure:perf_compare_with_existing_discrete}.
Our algorithm is able to solve 54 out of 80 benchmark instance to optimality, utilizing an average computation time of 274 seconds; of these, 22 instances were unsolved by the best previous method, while the remaining unsolved instances were determined to be within 1.2\% of optimality, on average.
As in the continuous setting, our algorithm strongly outperforms the existing method: it solves more instances and achieves the fastest computation time in all of them.

% Table generated by Excel2LaTeX from sheet 'compare - discrete'
\begin{table}[!htb]
  \centering
  \begin{threeparttable}
    \caption{Computational comparison of the proposed algorithm against the existing state-of-the-art algorithm (SD15)~\cite{Spliet2015_discrete} on all 80 benchmark instances of the discrete TWAVRP.}
    \begin{tabularx}{\textwidth}{lc*{6}{R}}
    \toprule
           &        & \multicolumn{3}{c}{SD15} & \multicolumn{3}{c}{Proposed algorithm} \\
    \cmidrule(lr){3-5}\cmidrule(l){6-8}
    $n$    & \#     & Optimal & Time (sec) & Gap (\%) & Optimal & Time (sec) & Gap (\%) \\
    \midrule
    10     & 10     & 10     & 3.9    & --               & 10     & 0.1    & -- \\
    15     & 10     & 10     & 185.9  & --               & 10     & 15.2   & -- \\
    20     & 10     & 9      & 1247.6 & 0.06             & 10     & 33.8   & -- \\
    25     & 10     & 3      & 504.4  & n/a$^\dagger$    & 9      & 248.8  & 0.43 \\
    30     & 10     & 0      & --     & n/a$^\dagger$    & 9      & 581.7  & 0.13 \\
    40     & 10     & 0      & --     & n/a$^\dagger$    & 5      & 1263.5 & 1.22 \\
    50     & 10     & 0      & --     & n/a$^\dagger$    & 1      & 533.6  & 1.31 \\
    60     & 10     & 0      & --     & n/a$^\dagger$    & 0      & --     & 1.30 \\
    All    & 80     & 32     & 457.5  & n/a$^\dagger$    & 54     & 274.4  & 1.21 \\
    Processor &        & \multicolumn{3}{r}{Intel Core i5-2450M 2.5 GHz} & \multicolumn{3}{r}{Xeon E5-2687W 3.1GHz} \\
    \bottomrule
    \end{tabularx}%
    \begin{tablenotes}%
    \small \item \textit{Note.} The reported results for ``SD15'' are the best entries of Tables~3 and~4 from that publication~\cite{Spliet2015_discrete}.
    \small \item $^\dagger$ The optimality gaps for the unsolved instances have not been reported in the publication.
    \end{tablenotes}%
  \label{table:computational_comparison_discrete}%
  \end{threeparttable}%
\end{table}%

\subsection{Detailed Discussion of Results}\label{sec:computational_results_detailed}
A comparison of Tables~\ref{table:computational_comparison_continuous} and~\ref{table:computational_comparison_discrete} shows that the discrete TWAVRP instances take longer to solve than the continuous ones.
This can be partly explained by the fact that, in the continuous setting, the separation problem~\eqref{eq:feasibility_milp} is a linear program, while in the discrete setting, it is a mixed-integer linear program.
Consequently, the algorithm spends a greater fraction of the total time in solving the separation problem in the latter case (see Table~\ref{table:percentage_of_computational_time}).

% Table generated by Excel2LaTeX from sheet 'detailed'
\begin{table}[!htb]
  \centering
  \caption{Percentage of computing time spent in various parts of the algorithm (averaged across instances solved to optimality in one hour). ``Solving VRPTW'' and ``Separation problem'' refer to Steps~\ref{path_based_algo:node_processing} and~\ref{main_algo:check_feasibility} of the algorithm from Section~\ref{sec:solution_approach_new_disjunctions} respectively, while ``Upper bounding'' refers to the steps in Section~\ref{sec:solution_approach_upper_bounds}.}
    \begin{tabularx}{\textwidth}{l*{6}{R}}
    \toprule
           & \multicolumn{3}{c}{Continuous} & \multicolumn{3}{c}{Discrete}  \\
    \cmidrule(lr){2-4}\cmidrule(l){5-7}
    $n$    & Solving VRPTW & Separation problem & Upper bounding & Solving VRPTW & Separation problem & Upper bounding \\
    \midrule
%    20        & 79.0   & 0.0    & 20.5   & 63.4   & 28.8   & 7.0 \\
    25        & 84.6   & 0.0    & 15.0   & 81.0   & 15.4   & 3.3 \\
    30        & 86.8   & 0.1    & 12.9   & 75.7   & 21.4   & 2.7 \\
    40        & 89.2   & 0.2    & 10.4   & 90.6   & 8.2    & 1.2 \\
    50        & 93.1   & 0.0    & 6.8    & 95.7   & 1.6    & 2.8 \\
    $\geq 25$ & 89.0   & 0.1    & 10.8   & 81.6   & 15.6   & 2.6 \\[0.1cm]
    \bottomrule    
    \end{tabularx}%
  \label{table:percentage_of_computational_time}%
\end{table}%

\begin{figure}[!htb]
\captionsetup[subfigure]{belowskip=0pt}
\centering
\begin{subfigure}[b]{0.5\textwidth}\centering
\includegraphics[height=6.2cm]{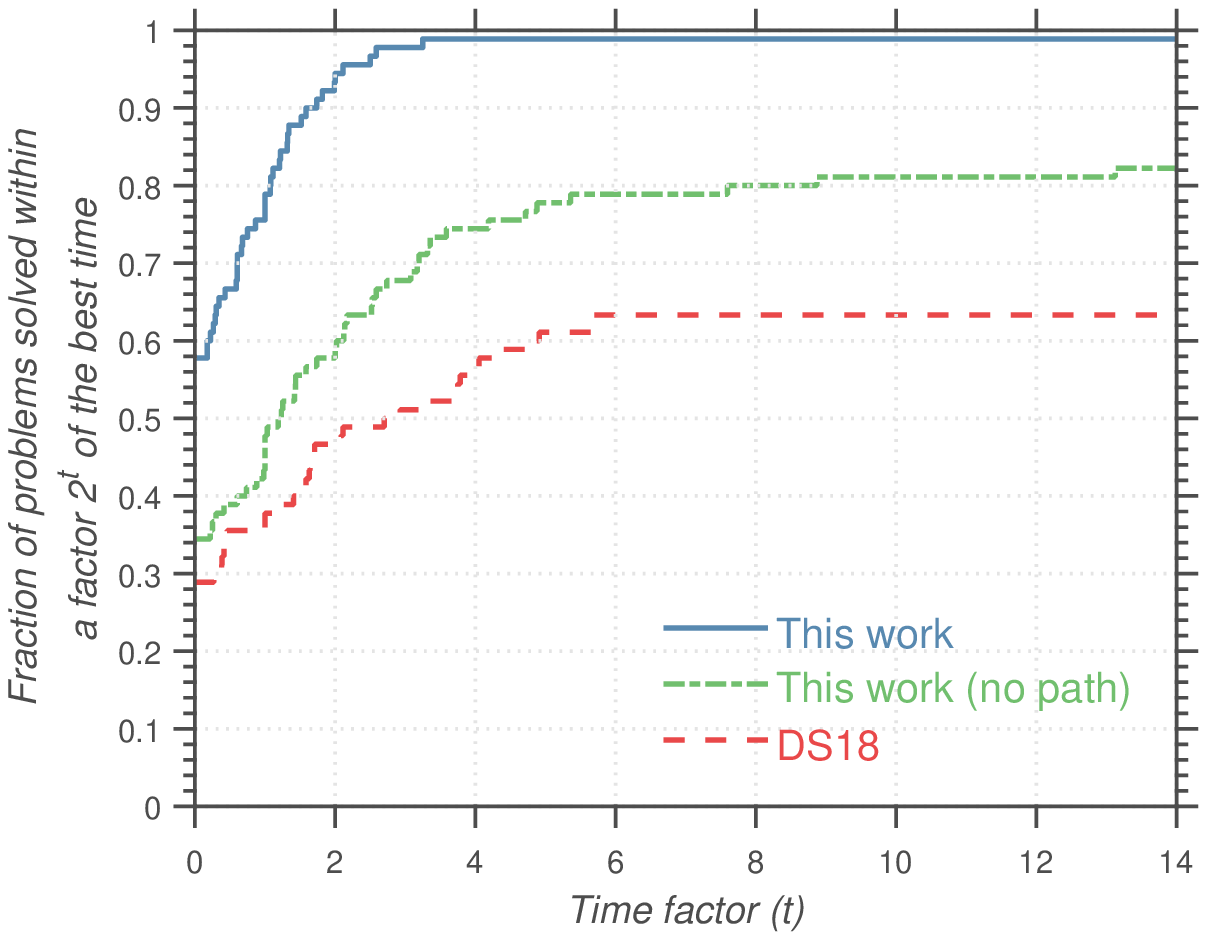}
\caption{Continuous}
\label{figure:perf_compare_with_tw_continuous}
\end{subfigure}~%
\begin{subfigure}[b]{0.5\textwidth}\centering
\includegraphics[height=6.2cm]{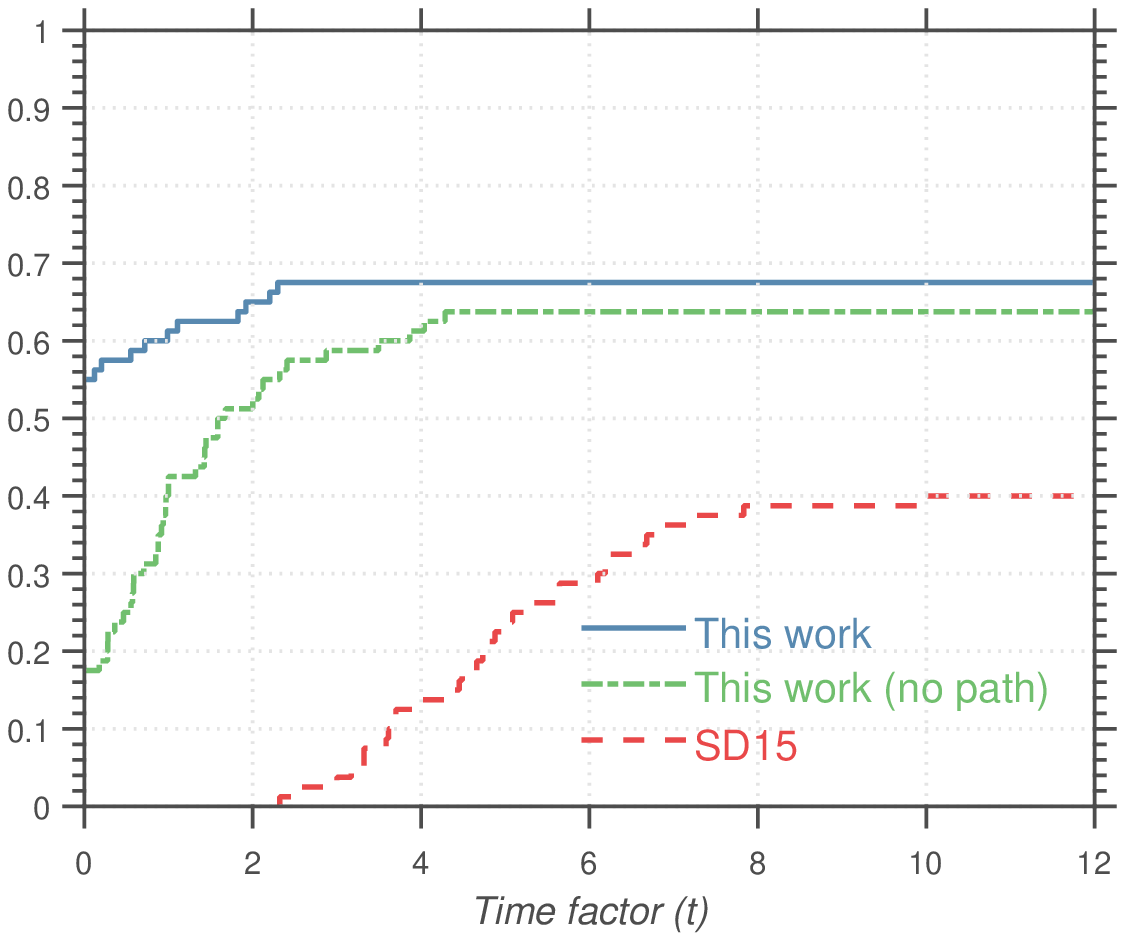}
\caption{Discrete}
\label{figure:perf_compare_with_tw_discrete}
\end{subfigure}
\caption{Log-scaled performance profiles across all benchmark instances. The profile ``This work (no path)'' refers to our proposed algorithm without the path-based disjunctions (see Section~\ref{sec:solution_approach_new_disjunctions}). The other profiles as well as the axes have the same meaning as in Figure~\ref{figure:perf_compare_with_existing}.}
\label{figure:perf_compare_with_tw}
\end{figure}

To show the efficacy of the path-based disjunctions and the associated branching rules (see Section~\ref{sec:solution_approach_new_disjunctions}), we disable them and run only the basic version of the algorithm from Section~\ref{sec:solution_approach_algorithm}. Table~\ref{table:computational_comparison_path_vs_tw} compares the performance of this basic version with the one incorporating the path-based disjunctions. They are also compared in the performance profiles of Figure~\ref{figure:perf_compare_with_tw}. The results indicate that the path-based branching rules are important to improve the tractability of the overall algorithm.
In particular, they are essential in reducing the total number of nodes that are processed (\textit{i.e.}, the number of times Steps~\ref{main_algo:node_processing} and~\ref{path_based_algo:node_processing} are executed in the overall search trees) from $2,400$ to $100$, on average.
We remark, however, that this reduction comes at a price: it requires modifying the underlying VRPTW solver (see Section~\ref{sec:implementation_details_vrptw_spp}). Nevertheless, even without the path-based branching rules, the basic version of our algorithm outperforms the existing ones (see Figure~\ref{figure:perf_compare_with_tw}), while having the advantage of being able to utilize any VRPTW solver in a modular fashion.

% Table generated by Excel2LaTeX from sheet 'compare - tw vs path'
\begin{table}[!htb]
  \centering
  \caption{Computational comparison of the algorithm with and without the path-based disjunctions on all 170 benchmark instances of the continuous and discrete TWAVRP.}
    \begin{tabularx}{\textwidth}{lc*{8}{R}}
    \toprule
           &        & \multicolumn{4}{c}{Without path-based disjunctions} & \multicolumn{4}{c}{With path-based disjunctions} \\
    \cmidrule(lr){3-6}\cmidrule(l){7-10}
    $n$                  & \#     & Optimal & Nodes  & Time (sec) & Gap (\%) & Optimal & Nodes  & Time (sec) & Gap (\%) \\
    \midrule
    $\left[10,15\right]$ & 40     & 40     &         5,887  &       91.5    & --     & 40     & 26     & 4.0    & -- \\
    $\left[20,25\right]$ & 40     & 36     &           933  &       147.7   & 0.17   & 39     & 158    & 68.7   & 0.43 \\
    $\left[30,35\right]$ & 30     & 22     &           949  &       327.6   & 0.16   & 28     & 132    & 220.7  & 0.46 \\
    $\left[40,45\right]$ & 30     & 19     &           469  &       210.7   & 0.66   & 25     & 139    & 534.1  & 1.22 \\
    $\left[50,60\right]$ & 30     & 8      &           557  &       1,061.0 & 1.12   & 11     & 105    & 681.2  & 1.30 \\[0.1cm]
    All                  & 170    & 125    &         2,427  &       229.4   & 0.75   & 143    & 108    & 208.9  & 1.19 \\
    \bottomrule
    \end{tabularx}%
  \label{table:computational_comparison_path_vs_tw}%
\end{table}%

Tables~\ref{table:detailed_performance_continuous} and~\ref{table:detailed_performance_discrete} present detailed results on all benchmark instances of the continuous and discrete TWAVRP, respectively. In these tables, if an instance could be solved to optimality within one hour, then \textbf{Opt [UB]} reports the corresponding optimal objective value, while \textbf{Time (sec) [LB]} reports the time to solve the instance to optimality. Otherwise, the columns respectively report in brackets the best upper and lower bounds found within the time limit of one hour.

% Table generated by Excel2LaTeX from sheet 'optimal values'
{\renewcommand{\arraystretch}{0.65}%
\begin{sidewaystable}
  \centering
  \caption{Detailed performance of the proposed algorithm on all 90 benchmark instances of the continuous TWAVRP.}
    \begin{tabularx}{\textwidth}{*{3}{ccCC}}
    \toprule
    Instance & $n$    & Opt\qquad [UB] & Time (sec) [LB] & Instance & $n$    & Opt\qquad [UB] & Time (sec) [LB] & Instance & $n$    & Opt\qquad [UB] & Time (sec) [LB] \\
    \cmidrule(r){1-4}\cmidrule(lr){5-8}\cmidrule(l){9-12}
    1      & 10     & 17.65  & 0.0    & 31     & 25     & 31.43  & 2.6    & 61     & 40     & 46.13* & 9.7 \\
    2      & 10     & 15.56  & 0.2    & 32     & 25     & 30.71  & 3.0    & 62     & 40     & 48.35  & 27.8 \\
    3      & 10     & 17.42  & 0.0    & 33     & 25     & 33.71  & 9.9    & 63     & 40     & 44.48* & 18.2 \\
    4      & 10     & 18.51  & 0.2    & 34     & 25     & 33.34  & 3.4    & 64     & 40     & 43.75  & 192.2 \\
    5      & 10     & 16.07  & 0.1    & 35     & 25     & 29.05  & 3.1    & 65     & 40     & 43.39* & 20.2 \\
    6      & 10     & 18.00  & 0.0    & 36     & 25     & 30.50  & 3.9    & 66     & 40     & 44.68* & 62.7 \\
    7      & 10     & 17.02  & 0.0    & 37     & 25     & 28.68  & 17.2   & 67     & 40     & 46.88* & 110.2 \\
    8      & 10     & 23.89  & 0.1    & 38     & 25     & 35.69  & 38.5   & 68     & 40     & 44.96* & 2,506.0 \\
    9      & 10     & 20.31  & 0.1    & 39     & 25     & 32.55  & 2.5    & 69     & 40     & 43.07* & 31.1 \\
    10     & 10     & 16.31  & 0.1    & 40     & 25     & 32.14  & 2.0    & 70     & 40     & 43.00* & 445.3 \\
    11     & 15     & 17.78  & 0.6    & 41     & 30     & 36.38  & 2.7    & 71     & 45     & 50.65* & 20.4 \\
    12     & 15     & 27.10  & 1.3    & 42     & 30     & 34.69* & 9.7    & 72     & 45     & 51.74* & 69.4 \\
    13     & 15     & 29.37  & 0.4    & 43     & 30     & 35.48  & 285.4  & 73     & 45     & 41.70* & 39.6 \\
    14     & 15     & 23.18  & 1.9    & 44     & 30     & 35.88  & 19.9   & 74     & 45     & 47.77* & 228.6 \\
    15     & 15     & 24.15  & 0.3    & 45     & 30     & 35.55  & 26.1   & 75     & 45     & 49.39* & 582.0 \\
    16     & 15     & 21.03  & 0.3    & 46     & 30     & 37.47  & 11.8   & 76     & 45     & 49.83* & 1,082.5 \\
    17     & 15     & 22.04  & 0.3    & 47     & 30     & 32.54  & 5.7    & 77     & 45     & 51.09* & 1,241.8 \\
    18     & 15     & 22.30  & 0.4    & 48     & 30     & 36.32  & 7.0    & 78     & 45     & 53.33* & 102.2 \\
    19     & 15     & 26.52  & 0.4    & 49     & 30     & 35.30  & 67.6   & 79     & 45     & 48.09* & 99.4 \\
    20     & 15     & 22.11  & 0.3    & 50     & 30     & 40.27  & 46.2   & 80     & 45     & 50.26* & 146.8 \\
    21     & 20     & 28.08  & 0.8    & 51     & 35     & 43.46  & 102.9  & 81     & 50     & 58.11* & 559.8 \\
    22     & 20     & 29.80  & 0.7    & 52     & 35     & 41.84  & 25.5   & 82     & 50     & 52.61* & 211.0 \\
    23     & 20     & 30.30  & 1.0    & 53     & 35     & 45.03* & 39.3   & 83     & 50     & 58.58* & 2,826.9 \\
    24     & 20     & 24.16  & 1.3    & 54     & 35     & 41.54* & 42.8   & 84     & 50     & 53.92* & 115.7 \\
    25     & 20     & 29.84  & 7.8    & 55     & 35     & 37.92  & 12.7   & 85     & 50     & 54.96* & 1,113.0 \\
    26     & 20     & 29.72  & 0.8    & 56     & 35     & 44.49* & 17.9   & 86     & 50     & 52.83* & 306.1 \\
    27     & 20     & 26.48  & 1.0    & 57     & 35     & [41.04]& [40.72]& 87     & 50     & 53.71* & 93.3 \\
    28     & 20     & 26.14  & 0.8    & 58     & 35     & 41.22  & 64.1   & 88     & 50     & 56.12* & 203.3 \\
    29     & 20     & 26.61  & 0.6    & 59     & 35     & 43.43  & 14.8   & 89     & 50     & 60.23* & 1,299.0 \\
    30     & 20     & 26.36  & 0.6    & 60     & 35     & 42.27  & 143.1  & 90     & 50     & 58.93* & 231.6 \\
    \bottomrule
    \multicolumn{12}{l}{\small *Instances solved to optimality for the first time are indicated with an asterisk.}
    \end{tabularx}%
  \label{table:detailed_performance_continuous}%
\end{sidewaystable}%
}

% Table generated by Excel2LaTeX from sheet 'optimal values'
{\renewcommand{\arraystretch}{0.65}%
\begin{sidewaystable}
  \centering
  \caption{Detailed performance of the proposed algorithm on all 80 benchmark instances of the discrete TWAVRP.}
    \begin{tabularx}{\textwidth}{*{3}{ccCC}}
    \toprule
    Instance & $n$    & Opt\qquad [UB] & Time (sec) [LB] & Instance & $n$    & Opt\qquad [UB] & Time (sec) [LB] & Instance & $n$    & Opt\qquad [UB] & Time (sec) [LB] \\
    \cmidrule(r){1-4}\cmidrule(lr){5-8}\cmidrule(l){9-12}
    1      & 10     & 12.83  & 0.1    & 31     & 25     & 35.47  & 24.3   & 61     & 50     & [52.23]& [51.57] \\
    2      & 10     & 16.84  & 0.2    & 32     & 25     & 32.66* & 16.0   & 62     & 50     & [55.61]& [55.19] \\
    3      & 10     & 16.60  & 0.1    & 33     & 25     & [31.75]& [31.61]& 63     & 50     & [50.75]& [50.09] \\
    4      & 10     & 15.96  & 0.2    & 34     & 25     & 34.14* & 236.9  & 64     & 50     & 51.17* & 533.6 \\
    5      & 10     & 19.65  & 0.2    & 35     & 25     & 30.29* & 617.7  & 65     & 50     & [54.11]& [53.44] \\
    6      & 10     & 18.13  & 0.3    & 36     & 25     & 32.54* & 611.6  & 66     & 50     & [57.52]& [56.57] \\
    7      & 10     & 12.17  & 0.1    & 37     & 25     & 27.48* & 568.5  & 67     & 50     & [58.14]& [57.54] \\
    8      & 10     & 17.09  & 0.2    & 38     & 25     & 34.83  & 15.9   & 68     & 50     & [56.37]& [55.27] \\
    9      & 10     & 20.14  & 0.1    & 39     & 25     & 34.39* & 123.7  & 69     & 50     & [53.85]& [53.40] \\
    10     & 10     & 17.17  & 0.1    & 40     & 25     & 30.73  & 25.0   & 70     & 50     & [57.37]& [56.37] \\
    11     & 15     & 23.04  & 7.4    & 41     & 30     & 36.39* & 37.9   & 71     & 60     & [64.83]& [63.60] \\
    12     & 15     & 25.27  & 1.0    & 42     & 30     & 40.59* & 245.4  & 72     & 60     & [62.60]& [61.57] \\
    13     & 15     & 22.12  & 2.8    & 43     & 30     & 37.18* & 88.7   & 73     & 60     & [64.92]& [63.91] \\
    14     & 15     & 18.46  & 0.7    & 44     & 30     & [38.02]& [37.97]& 74     & 60     & [69.14]& [68.59] \\
    15     & 15     & 24.87  & 129.9  & 45     & 30     & 36.72* & 311.2  & 75     & 60     & [63.61]& [63.15] \\
    16     & 15     & 19.82  & 2.4    & 46     & 30     & 34.76* & 237.8  & 76     & 60     & [64.49]& [63.77] \\
    17     & 15     & 21.96  & 4.7    & 47     & 30     & 42.24* & 133.6  & 77     & 60     & [61.24]& [60.82] \\
    18     & 15     & 22.93  & 0.7    & 48     & 30     & 37.04* & 3,501.3& 78     & 60     & [64.77]& [63.06] \\
    19     & 15     & 23.14  & 1.3    & 49     & 30     & 40.47* & 202.6  & 79     & 60     & [65.48]& [64.91] \\
    20     & 15     & 18.84  & 0.7    & 50     & 30     & 39.89* & 477.2  & 80     & 60     & [64.42]& [63.76] \\
    21     & 20     & 27.99  & 1.2    & 51     & 40     & [41.96]& [41.44]&        &        &        &  \\
    22     & 20     & 25.63  & 16.8   & 52     & 40     & [47.43]& [47.37]&        &        &        &  \\
    23     & 20     & 26.53  & 199.2  & 53     & 40     & 41.76* & 1,829.8&        &        &        &  \\
    24     & 20     & 32.36  & 3.9    & 54     & 40     & 45.96* & 1,385.2&        &        &        &  \\
    25     & 20     & 28.84  & 4.7    & 55     & 40     & [48.68]& [48.10]&        &        &        &  \\
    26     & 20     & 26.99  & 6.3    & 56     & 40     & [44.88]& [43.99]&        &        &        &  \\
    27     & 20     & 27.55* & 80.0   & 57     & 40     & 43.90* & 619.3  &        &        &        &  \\
    28     & 20     & 26.53  & 14.4   & 58     & 40     & 43.09* & 918.1  &        &        &        &  \\
    29     & 20     & 29.49  & 8.6    & 59     & 40     & [49.27]& [48.50]&        &        &        &  \\
    30     & 20     & 23.55  & 2.5    & 60     & 40     & 47.13* & 1,565.2&        &        &        &  \\
    \bottomrule
    \multicolumn{12}{l}{\small *Instances solved to optimality for the first time are highlighted with an asterisk.}
    \end{tabularx}%
  \label{table:detailed_performance_discrete}%
\end{sidewaystable}%
}

\subsection{Instances containing a Large Number of Scenarios}\label{sec:computational_results_large_instances}
We now turn our attention to benchmark instances containing a large number of scenarios.
Our goals are two-fold.
First, we aim to understand how our algorithm performs as the number of considered scenarios ($S$) increases.
Second, we aim to understand the cost benefits of considering more scenarios during strategic time window assignment.
In pursuit of these goals, we consider the 80 benchmark instances consisting of 15 demand scenarios each (see Section~\ref{sec:computational_results_instances}).
For each of these instances, we obtain time window assignments using our algorithm by considering the following sample average approximations: \textit{(i)} the original instance with all 15 scenarios, and \textit{(ii)} the original instance with only the first $S$ scenarios, where $S \in \{1, 3, 5, 10\}$.
For each approximation, we implement a parallel version of our algorithm in which Step~\ref{path_based_algo:node_processing} and the upper bounding step in Section~\ref{sec:solution_approach_upper_bounds} are each parallelized using up to $\min\{S, 10\}$ threads.

\begin{figure}[!htb]
\centering
\caption{The ratio of CPU time to wall clock time with increasing number of scenarios (for instances whose solution took at least one wall clock second and at most one wall clock hour). For each $S$, the lower and upper most dashes denote the minimum and maximum values of the ratio, the lower and upper edges of the box denote the first and third quartile, while the middle bar (in red color) denotes the median ratio.}
\label{figure:parallel_speedup}
\includegraphics[scale=0.75]{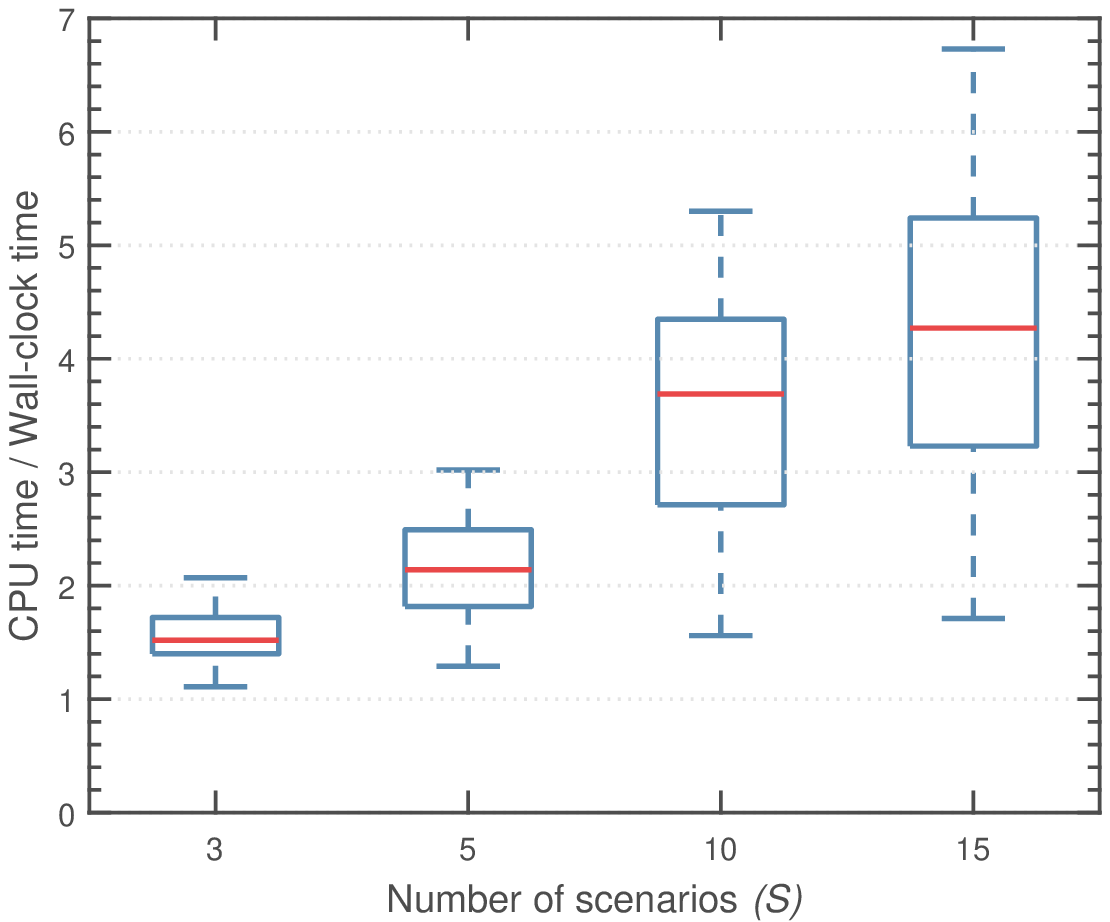}
\end{figure}

For $S=1$ and $S=3$, all 80 instances were solved to optimality, while 73 out of 80 instances were solved to optimality for $S = 5$ (the average gap for the 7 unsolved instances being less than 0.2\%). Therefore, at the interest of brevity, we shall not show tabulated results for these cases. On the other hand, for the cases with the higher number of scenarios, namely $S=10$ and $S=15$, Table~\ref{table:performance_large_S} presents a summary of the performance. The columns in this table have the same interpretation as in Tables~\ref{table:computational_comparison_continuous} and~\ref{table:computational_comparison_discrete}; the only difference is that the column \textbf{Time (sec)} is now broken into two parts: \textbf{Wall} denotes the average wall clock time, while \textbf{CPU} denotes the average CPU time.
The ratio between these quantities is a good measure of how well the algorithm scales across multiple threads (i.e., how much it benefits from parallelism), which is also plotted as a function of $S$ in Figure~\ref{figure:parallel_speedup}.

% Table generated by Excel2LaTeX from sheet 'correlated - table'
\begin{table}[!htb]
  \centering
  \caption{Summary of computational performance of the parallelized algorithm on 80 benchmark instances of the (continuous and discrete) TWAVRP, each containing $S$ scenarios.}
    \begin{tabularx}{\textwidth}{lc*{8}{R}}
    \toprule
           &        & \multicolumn{4}{c}{$S = 10$}          & \multicolumn{4}{c}{$S = 15$} \\
    \cmidrule(r){3-6}\cmidrule(l){7-10}
           &        &        & \multicolumn{2}{c}{Time (sec)} &  &  &  \multicolumn{2}{c}{Time (sec)} & \\
    \cmidrule(lr){4-5}\cmidrule(lr){8-9}
    $n$    & \#     & Optimal& Wall   & CPU    & Gap \% & Optimal& Wall    & CPU     &  Gap \% \\
    \midrule
    10     & 20     & 20     & 0.7    & 3.8    & --     & 20     & 1.1     & 7.0     & -- \\
    15     & 20     & 18     & 127.7  & 474.6  & 0.21   & 18     & 441.7   & 1,639.2 & 0.64 \\
    20     & 20     & 16     & 266.4  & 779.8  & 0.46   & 11     & 426.2   & 1,403.3 & 0.33 \\
    25     & 20     & 5      & 465.3  & 841.6  & 0.60   & 2      & 2,189.9 & 3,777.1 & 0.70 \\[0.1cm]
    All    & 80     & 59     & 150.9  & 428.9  & 0.53   & 51     & 334.2   & 1,032.1 & 0.58 \\
    \bottomrule
    \end{tabularx}%
  \label{table:performance_large_S}%
\end{table}%

Table~\ref{table:performance_large_S} shows that we can consistently solve all benchmark instances with 15 scenarios to an optimality gap of less than 1\% within a time limit of one hour. Figure~\ref{figure:parallel_speedup} shows that the speedup in wall clock time is sublinear with respect to the number of parallel threads. The reasons for deviating from a perfect linear speedup are two-fold. First, each node of our search tree does not necessarily require the solution of $S$ VRPTW instances (see concluding paragraph of Section~\ref{sec:solution_approach_algorithm}).
Indeed, for the considered benchmark instances, the average number of VRPTW instances solved in a typical node is smaller than $S/2$.
Second, the variance in solution times across the VRPTW instances solved in a node is typically large because the feasible route sets in a particular scenario might be drastically different compared to other scenarios in the same node (because of different demand realizations).
Consequently, the different CPU threads are not necessarily balanced, i.e., do not perform equal amount of work.
Nevertheless, as Figure~\ref{figure:parallel_speedup} shows, on average, our algorithm performs four times as many computations for instances containing 15 scenarios by utilizing up to ten threads as compared to using just one.

Finally, Table~\ref{table:simlation_results} shows the cost savings in routing that are to be expected by considering more than just one scenario during time window assignment.
We calculate these savings as follows.
First, for a particular instance (with a given $S$), we obtain the optimal time window assignment $\tau^\star$ from our algorithm.
Next, we calculate the ``in-sample'' expected costs by \textit{(i)} fixing the time windows to $\tau^\star$, \textit{(ii)} solving a VRPTW instance with time windows fixed to $\tau^\star$ for each of the 15 postulated demand scenarios in the original benchmark instance, and \textit{(iii)} averaging the costs.
The ``out-of-sample'' expected costs are obtained in exactly the same manner except that in step \textit{(ii)}, the costs are evaluated over 100 independently generated demand scenarios, which are randomly drawn using the procedure described in Section~\ref{sec:computational_results_instances}.
Table~\ref{table:simlation_results} shows that, on average, we expect to do better than the deterministic solution by about $2.3\%$ when considering up to $3$ scenarios and by about $3.2\%$ when considering up to $15$ scenarios (based on out-of-sample evaluations).
If we consider only the discrete instances, the out-of-sample savings increase to about $3.0\%$ for $S=3$ and $3.7\%$ for $S=15$.
In either case, we observe that the marginal benefits diminish as the number of scenarios grows.
Although we expect this trend to generally hold, we also expect the actual magnitudes of the cost savings to depend on the problem parameters as well as the dimensionality of the uncertainty (the number of customers in this case).
Indeed, for high-dimensional problems with narrow and discrete time windows, we expect greater cost benefits by considering more scenarios of the uncertainty.

% Table generated by Excel2LaTeX from sheet 'correlated - eval'
\begin{table}[!htb]
  \centering
  \caption{Expected cost savings from considering $S$ scenarios during strategic time window assignment relative to considering only one scenario (averaged across all 60 instances with $n \geq 15$).}
    \begin{tabularx}{0.8\textwidth}{l*{5}{C}}
    \toprule
                  & $S=1$  & $S=3$  & $S=5$  & $S=10$ & $S=15$ \\
    \midrule
    In-sample     & 0.00\% & 2.26\% & 2.88\% & 3.13\% & 3.22\% \\
    Out-of-sample & 0.00\% & 2.36\% & 2.89\% & 3.18\% & 3.21\% \\
    \bottomrule
    \end{tabularx}%
  \label{table:simlation_results}%
\end{table}%

\section{Conclusions}\label{sec:conclusions}
While there has been significant research done on efficient algorithms for vehicle routing problems, rapidly changing business models and advances in technology are creating new challenges that are yet unsolved.
This paper attempts to address the challenge of dealing with operational uncertainty during strategic decision-making in the context of vehicle routing.
In particular, we studied problems in which the decisions correspond to an allocation of long-term delivery time windows to customers.
These problems are motivated from several real-world distribution operations, and are particularly common in retail.
We proposed a novel algorithm to solve this problem that is highly competitive with existing methods. It draws on existing algorithms for deterministic vehicle routing problems (both exact and heuristic) and can use any vehicle routing solver as a subroutine, thus facilitating its deployment in practice.
From a modeling viewpoint, it allows the user to postulate potential scenarios of future uncertainty corresponding to different routing-specific parameters, as well as incorporate modular changes in routing-specific constraints.
Our business insights, aided via numerical experiments, are that long-term costs are expected to decrease by modeling more scenarios of future uncertainty but the marginal benefits rapidly diminish as a function of the number of postulated scenarios. In other words, a few scenarios are sufficient to obtain time window assignments that would incur lower long-term costs compared to assignments that would be obtained by completely ignoring operational uncertainty.

Our results open multiple avenues for future research. First, the scope of the proposed algorithm can be broadened further by incorporating features such as time-dependent travel times (which capture predictable variations in traffic conditions) as well as soft time windows (which allow violation of the assigned time windows in a controlled manner). Second, the time window assignment model can be generalized to allow customers to express preferences in the time window allocation process. More ambitiously, it would be instructive to explore how we can reinterpret these concepts in the context of consistent vehicle routing problems.

\section*{Acknowledgments}
We acknowledge support from the United States National Science Foundation, award number CMMI-1434682.
Anirudh Subramanyam and Akang Wang also acknowledge support from, respectively, the John and Claire Bertucci Graduate Fellowship program and the James C. Meade Fellowship program, at Carnegie Mellon University.

\bibliographystyle{plain}
\bibliography{bibliography}
\end{document}